\documentclass[12pt]{amsart}

\usepackage{amssymb,amsmath}
 \usepackage[dvips]{graphicx}
\usepackage{verbatim}

\makeatletter
\def\@cite#1#2{{\m@th\upshape\bfseries%
[{#1\if@tempswa{\m@th\upshape\mdseries, #2}\fi}]}}
\makeatother
\theoremstyle{plain}
\newtheorem{thm}{Theorem}[section]
\newtheorem{lma}[thm]{Lemma}
\newtheorem{cor}[thm]{Corollary}
\newtheorem{prop}[thm]{Proposition}
\theoremstyle{definition}

\newtheorem{defn}[thm]{Definition}

\newcommand{\Q}{{\mathbb{Q}}}
\newcommand{\upchi}{{\raise.35ex\hbox{$\chi$}}}
\begin{document}

\title[Generic $3$-connected planar constraint systems]%
{Generic $3$-connected planar constraint systems are not soluble by radicals}
%

\author[J.C. Owen]{J.C. Owen}
\address{D-Cubed Ltd\\
         Park House\\
         Cambridge CB3 0DU\\
        United Kingdom}
\email{john.owen@d-cubed.co.uk}
\author[S.C. Power]{S.C. Power}
\address{Department of Mathematics and Statistics\\
         Lancaster University\\
         Lancaster, LA1 4YF\\
         United Kingdom}
\email{s.power@lancaster.ac.uk}

\thanks{{\it key words and phrases. maximally independent graph, 3-connected,
algorithms for CAD, solvable by radicals.} }

\date{}

\begin{abstract}
We show that planar embeddable $3$-connected CAD 
graphs are generically non-soluble.
A CAD graph represents a configuration of points on the Euclidean plane with just
enough distance dimensions between them to ensure rigidity. 
Formally, a CAD graph is  
a maximally independent graph, that is, one that satisfies
the  vertex-edge count  $2v - 3 = e$ together with a
corresponding inequality for each subgraph. The following
main theorem of the paper resolves 
a conjecture of Owen \cite{Owen} in the planar case.
Let $G$ be a maximally independent $3$-connected planar
graph, with more than 3 vertices, together with a realisable
assignment of generic dimensions for the edges which includes 
a normalised unit length (base) edge.
Then, for any
solution configuration for these dimensions on a plane, with the
base edge vertices placed at rational points,  not all
coordinates of the vertices lie in a radical extension of the
dimension field.
\end{abstract}

\maketitle
\setcounter{tocdepth}{1}
\tableofcontents

\section{Introduction}

A fundamental problem in Computer Aided Design (CAD) is the
formulation of effective algebraic algorithms or numerical
approximation schemes  which solve for the location of points on a
plane, given a set of relative distances between
them. 
The relative distances are
usually called dimensions in CAD by analogy with the dimensions on a
dimensioned drawing and we will adopt that terminology here.
For CAD applications the relevant class of configurations
are those for which the dimensions are just sufficient to ensure
that the points are located rigidly with respect to one another.
It is a well known result of Laman \cite{lam} that the graphs
underlying generically rigid configurations (frameworks) have a
simple combinatorial description.  In our terminology they are the
 so-called maximally independent graphs, that is, those
satisfying the vertex-edge count  $2v - 3 = e$ together with a
corresponding inequality for each subgraph.

A number of algebraic and numerical methods have been proposed for
solving these plane configurations (Owen \cite{Owen}, Bouma et al
\cite{bouma-et-al}, Light and Gossard \cite{light-gossard}) and
these have been successfully implemented in CAD programs.
Algebraic and combinatorial algorithms for graphs are particularly
desirable for their speed and robustness and the resulting
dramatic efficiency gains. For instances of this 
see, for example, the quadratic
extension algorithm of \cite{Owen}, the graph decomposition
algorithm of Hopcroft and Tarjan \cite{hop-tar}, or  the
combinatorial approach to protein molecule flexibility in Jacobs
et al \cite{jac}.

 Current algebraic methods for solving CAD graphs
assemble the solution for complete configurations from the
solutions of rigid subcomponents and the assembly process involves
only rigid body transformations, fusion at vertex pairs, and  the
solution of quadratic equations. The simplest subcomponent is a
triangle of points which is solvable by quadratic equations and it
thus follows that if the original configuration is assembled from
triangles then it is solvable through
successive quadratic extensions of the
ground dimension field. The other subcomponents possible in this
process are all represented by graphs which are  $3$-connected (in
the usual sense of vertex $3$-connected
 \cite{Tutte}) and so the problem of solving
general configurations passes to the problem of solving
configurations which are represented by $3$-connected graphs.
Determination of $3$-connectivity can be effected rapidly with
order $O(v + e)$. (See \cite{hop-tar}.)

We have previously suggested that with generic dimension values a
subcomponent which is represented by a $3$-connected graph cannot
be solved by quadratic equations (Owen \cite{Owen}).
Configurations that can be solved in this way are also known as
"ruler and compass constructible" and Gao  and Chou
\cite{Gao-Chou} have given a procedure for determining in
principle if any given configuration is ruler and compass
constructible. However their analysis is based on the detail of
derived elimination equations and they do not address the problem
of generic solubility or non-solubility  for general classes of
graphs.

Despite the importance of algebraic solubility, the
intractability or otherwise of generic $3$-connected
configurations has not been put on a firm theoretical basis and in
the present paper we begin such a project.

The solution configurations that we consider 
are comprised of points in the plane with a number of
specified distances (dimensions) between them. With the natural correspondence
of points to vertices and constraint pairs to edges each constraint system
has an associated abstract graph. 
It is the nature of the abstract graph that is significant for the solubility
of the constraint system and we shall be concerned with
the situation where the abstract graph is a planar
graph in the usual graph-theoretic sense; it can be drawn with
edges realised by curves in the plane with no crossings.

 We show that a planar
$3$-connected maximally independent graph with generic
dimensions is not only not solvable by quadratic extensions but is
not soluble by radical extensions, that is, by means of the
extraction of roots of arbitrary order together with the basic
arithmetical operations. In fact our methods make use of some
intricate planar graph theory leading to a edge contraction
reduction scheme which is also of independent interest. The main
theorem of the paper can be stated as follows.
\medskip

\begin{thm}\label{T:main}
{\it Let $G$ be a maximally independent $3$-connected planar
graph, with more than 3 vertices, together with a realisable
assignment of generic dimensions for the edges which includes 
a normalised unit length (base) edge.
Then, for any
solution configuration for these dimensions on a plane, with the
base edge vertices placed at rational points,  not all
coordinates of the vertices lie in a radical extension of the
dimension field.}
\end{thm}

\medskip

It follows in particular that the current algebraic schemes already solve
all of the generic configurations with a planar graph that can be solved by
radical extensions !  Also, we conjecture that planarity is not
necessary for this conclusion.

Recall that a celebrated and fundamental achievement of classical
Galois theory is that a polynomial of degree 5 or more, with
rational coefficients, is not generally soluble by radical
extensions over ${\mathbb Q}$. For a generic version of this, one
can assert that a generic monic polynomial of degree $r \ge 5$ is
not soluble by radical extensions of the base field ${\mathbb
Q}(\{d\})$, where $\{d\} = \{d_1,...,d_{r-1}\}$ are the generic
(algebraically independent) coefficients. In this case, with
coefficient field understood, the polynomial is said to be,
simply, non-soluble. These facts suggest that if one is presented,
as we are here, with $N$ polynomial equations in $N$ unknowns,
with no apparent step by step solution scheme involving at most
degree 4 polynomials, then solutions will not lie in radical
extensions of the coefficient field. On the other hand, possibly
working against this intuition is the fact that our constraint
equations are all of quadratic type, in four variables, with a
single generic constant term, and the variables of the equations
reflect a (planar) graph structure which may possess an intrinsic reduction
scheme. However our result shows that in fact  there can be no grounds for
a solution scheme by radical extraction which embraces more than
the known quadratically soluble graphs. To paraphrase Theorem 1.1,
{\it  planar embeddable $3$-connected CAD graphs are generically non-soluble.}

Let us now outline the structure of the proof, the
entirety of which is lengthy and eclectic, making use of graph
theory, elimination theory for the ideals of complex affine
varieties, Galois theory for specialised coefficient fields, and a
brute force demonstration of the non-solubility of a 
vertex minimal
$3$-connected maximally independent 
planar graph. We refer to this graph, indicated in
Figure 1, as the doublet.

\begin{figure}
\centering
\includegraphics[width=9cm]{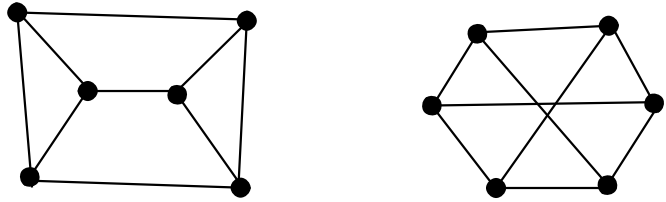} \caption{The Doublet and $K_{33}$.}
\end{figure}
\medskip

The fact that the generic doublet graph is not soluble by radicals 
is obtained in Section 8 by first obtaining 
an explicit integral dimensioned doublet
which is not soluble. Here the Galois groups of univariate
polynomials in the elimination ideals for the constraint equations are
computed with some computer algebra assistence.
Generic non-solubility then follows from our 
Galois group specialisation theorem.

The strategy of the proof is to show
that if there exists a graph $G$ which is maximally independent, planar, 
3-connected 
and radically soluble then there is a smaller such graph
with fewer vertices. By the minimality of
the doublet this implies that the doublet is radically soluble
which gives  the desired contradiction

There are two aspects  to the reduction step.
The first of these is purely graph theoretic
and is dealt with in the extensive analysis of Section 4.
The main theorem
there shows that a 3-connected planar
maximally independent graph $G$ has either an edge
$e$ in a triangle of edges which can be contracted to give
 a smaller such graph $G/e$, or has a rigid subgraph
which can be replaced by a triangle to produce a smaller such graph,
$H$ say. 
The second aspect is to connect the solubility of the (finite) variety
of solutions for the dimensioned graph $G$ to that
of the varieties of the resulting smaller dimensioned graphs.
In the latter case we can simply compare generic constraint equations
(see Proposition \ref{P:8.1})
to deduce that
\[
\mbox{generic } G \mbox{ radical } \Rightarrow \mbox{generic } H \mbox{ radical }
\]
However the former case of edge contraction is much more subtle.
We approach this by
 noting first that the complex variety $V(G/ e)$ of
solutions for the generic contracted graph is  identifiable with
the variety of solutions for $G$ with partially specialised
 dimensions, with the contracted edge dimension $d_e$
specialised to $0$ and the two other edges of the contracted
triangle  specified as being equal. 
This gives the easy implication
\[
\mbox{specialised } G \mbox{ radical } \Rightarrow \mbox{generic } G/e \mbox{ radical }
\]
However we now need the final step, that is the implication
\[
\mbox{generic  } G \mbox{ radical } \Rightarrow \mbox{specialised } G \mbox{ radical }
\]
To obtain this
we consider carefully
the polynomials which are the generators of the
single variable elimination ideals associated with  
the constraint equations. We relate these generators
to the corresponding polynomials for the ideals of the specialised
equations. In fact we relate the solubility of these polynomials
through a two-step process for the double specialisation. This is
effected in Sections \ref{S:elimination}, \ref{S:reduction}.
The proof of the final step is then completed by means of 
another application of the
Galois group specialisation theorem, Theorem \ref{T:galois2}. This
theorem asserts, roughly speaking,
 that the Galois group of a polynomial $p$ is a subgroup of
the Galois group of a polynomial $P$ when $p$ derives from $P$ by
partial specialisation of coefficients. We were unable to find a
reference for this seemingly classical assertion.

Let us highlight two very important
ideas which run through the
proof of the reduction step for edge contractions
(Theorem \ref{T:reduction}).

The first of these is that we must restrict attention to graphs
whose constraint equations, both generic and
specialised, have finitely many complex solutions.
This form of rigidity
for complex variables we call {\it zero dimensionality} and its 
significance is explained
fully in the next section. It guarantees that univariate
elimination ideals for the constraint equations are generated by
univariate polynomials. Unfortunately, to maintain zero dimensionality 
our contraction scheme to the doublet must  operate
entirely in the framework of maximally independent graphs and it is this that 
necessitates  the extended graph theory of Section 4.

The second important idea is that the constraint equations happen to be of parametric
type.
As is well known this means that various associated complex affine varieties
are irreducible and in particular (Theorem \ref{T:Vb}) this 
is so for the so-called big variety in which the coordinates of
vertices and the dimensions of edges are viewed as complex variables.
With irreducibility present we can arrange the univariate generators of single
variable elimination ideals to be irreducible over the appropriate field  
(Theorem \ref{T:generator}) and so
either all roots of the generator are radical or none are.
Now it is the case that not every root of the generator need 
derive from a solution of the
constraint equations. Thus the fact that either all roots or no roots are radical
allows us to compare the solubility or otherwise of
$G$ and $G/e$ by examining the solubility or otherwise of these univariate generators
(Theorem \ref{T:AG4a}).

Finally, we remark that the assumption that graphs have a planar
embedding is used to guarantee that there is a reduction scheme to
a minimal graph based on contracting edges. We expect that
there are more general reduction schemes which terminate in either
the doublet or the non-planar graph $K_{33}$. Also we are able to
show that $K_{33}$ is generically non-soluble and this gives
further support to our conjecture that general $3$-connected
maximally independent graphs are non-soluble.

The results of this article were announced at the Fourth
International Workshop on Automated Deduction in Geometry in
September 2002 \cite{owe-pow-1}. We thank Walter Whiteley for
helpful discussions and for directing our attention to the paper
of Asimov and Roth \cite{asi-rot}
\section{Constraint equations and algebraic
varieties}\label{S:constraint}

\medskip
We begin by formulating the main problem which is to determine
the complex algebraic variety arising from the solutions to
the constraint equations of a normalised dimensioned graph.

Let $G = (V,E)$ be a graph with vertex set $V$ and edge set $E$.
We are concerned with the problem of
determining coordinates $(x_v,y_v)$ for each vertex $v$ so that for
some preassigned
dimensions  for the edges $e$ in $E$, we have solutions
to the set of equations
\[ f_e = 0, ~~~~~~e \in E,
\]
where, for the edge $e = (vw)$,
\[
f_e = (x_v - x_w)^2 + (y_v - y_w)^2 - d_e.
\]
The
dimensions $d_e$ are taken to be nonnegative real numbers,
representing the square of the edge lengths of realised graphs.

It is convenient to refer to the set $\{f_e\}$ as a set of
(unnormalised) {\it constraint equations} for the graph.
Although in practice one is interested primarily in the real
solutions in $\mathbb{R}^2$ for the vertices, which in turn
account for the Euclidean realisations of the dimensioned graph,
it is essential to our approach that we consider all  complex
solutions. In this case solutions always exist and we can employ
the elimination theory for complex algebraic varieties.

Bearing in mind the multiplicity of solutions associated with
Euclidean isometries we assume that for some base edge $b = (vw)$
in $E$ we have $d_b = 1$ and the specification $(x_v,y_v) = (0,0),
(x_w,y_w) = (1,0).$ This gives rise to a set of normalised
constraint equations $\{f_e\}$. If, in addition, the dimensions are
algebraically independent then we say that $\{f_e\}$ is a set of
{\it generic constraint equations} for $G$.
We shall generally assume that dimension sets and sets of
constraint equations are normalised.

Let $(G,\{d_e\})$ be a normalised dimensioned graph with $n$
vertices and let $x_i, y_i, 1 \le i \le n-2$, be the coordinate
variables for the non-base vertices. We write $V(\{f_e\})$ for the
complex affine variety in $\mathbb{C}^{2n-4}$ determined by the
corresponding set of constraint equations $\{f_e\}$.
\medskip

We now give  some definitions which 
give precise meanings to the terms
generic and rigid. 
There is a close connection between our
formalism
and that of the theory of rigid frameworks (see Whiteley \cite{whi-2}
and  Asimow and Roth \cite{asi-rot})
and in particular the notion 
of an independent graph is taken from this context.

\begin{defn}\label{D:zerodim}
The dimensioned graph  $(G,\{d_e\})$ is said to be
{\it zero dimensional} if the complex
algebraic variety $V(\{f_e\})$ is zero dimensional, that is, $V(\{f_e\})$
is a finite non-empty set.
\end{defn}

\begin{defn}\label{D:indptG}
Let $G$ be a graph with $v_G$ vertices and $e_G$ edges. Then $G$
is said to be {\it independent} if for every vertex induced subgraph
$H$, we have $2v_H - e_H \ge 3$. The graph $G$ is said to 
be {\it maximally independent} if it is independent and in
addition $2v_G - e_G = 3$.
\end{defn}

The graphs for which generic dimensions give zero dimensional
varieties admit a simple combinatorial description as we see 
below. These maximally independent graphs are also
known colloquially as CAD graphs.
This equivalence 
follows from our variant of  Laman's theorem. In fact we shall only need one direction,
proved in Theorem \ref{T:Laman} namely that maximally independent graphs with
generic dimensions are zero dimensional.


Let us indicate  more fully the nature and  significance 
of zero dimensionality.

Zero dimensionality for dimensioned graphs might
also be termed {\it complex rigidity}. For {\it non-generic} dimensions
it is a stronger requirement than the rigidity of the graph as a
bar-joint structure as given in 
 \cite{asi-rot} and \cite{whi-2}.
To appreciate  this consider the maximally independent graph in
Figure 2 which we view as a generically dimensioned graph with
normalised dimensions $\{d_e\}$. The two arrowed edges suggest a
specialisation of $\{d_e\}$ to a new dimension set for the same
graph in which  the arrowed edges have length zero and  two pairs
of edges are of equal length. Despite the fact that the resulting
semi-generic bar-joint structure is physically rigid and that the
original graph has been contracted  onto a maximally independent
graph (the doublet), the specialised dimensioned graph is {\it not}
zero-dimensional. In this case the variety $V'$ is
a one-dimensional variety in $\mathbb{C}^{12}$
which  meets the real subset $\mathbb{R}^{12}$  in a finite set. 

The graph in Figure 2 is not $3$-connected.
In fact it  is quadratically soluble in the sense expressed in Theorem 3.2. 
However the doublet is
$3$-connected and, as we will show in a subsequent section, 
it is not quadratically
soluble. This observation indicates that in any reduction scheme
for the proof involving edge contractions it is necessary
to work within the
category of zero dimensional graphs rather than rigid graphs in
the usual sense.

\begin{figure}
\vspace{.5in}
\centering
\includegraphics[width=7cm]{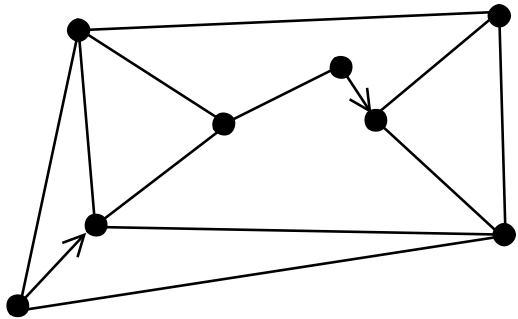}
\caption{Specialisation of a quadratically soluble graph onto a
doublet.}
\end{figure}

.

\pagebreak

The following general theorem will be used in the proof of Theorem \ref{T:Laman}.
By a specialisation of the dimension set $\{d_e\}$ in $\mathbb{C}^n$
(generally  an algebraically independent set) we mean a set
$\{d_e'\}$ in which some or all of the $d_e$ have been replaced by
rational numbers.

\begin{thm}\label{T:multispecialisation}
 Let $V$ be a complex affine variety in $\mathbb{C}^n$ defined by polynomial
equations of the form
\[
f_i = h_i(\{x_j\})-d_i = 0,~~~~ 1 \le i \le n,
\]
where $\{h_i\}$ are polynomials with rational coefficients
in the complex variables $\{x_j\}= \{x_1,...,x_n\}$, and where $\{d_i\}$ 
a set of constants in
$\mathbb{C}$. If $J$ is the $n \times n$ matrix 
$J=(J_{ij}) = (df_i/dx_j)$ and $\det(J)$ is not identically
zero as a polynomial in $\{x_j\}$ then

(1) The coordinates $\{x_j\}$ of any zero of $V$ are algebraically independent as a
set if and only if the constants $\{d_i\}$ are algebraically independent as a
set.

(2) If $\{d_i\}$ are algebraically independent then $\dim (V)=0$.
\end{thm}

\begin{proof}
Suppose that $\{d_i\}$ are algebraically dependent then
there is some polynomial $p$ in $n $ variables 
with  $p(d_1,\dots ,d_n)=0.$ Define the polynomial
$q$ by 
$q(\{x_j\})=p(h_1(\{x_j\}),\dots ,h_n(\{x_j\}))$. Then $q$ is not
the zero polynomial because $\det(dh_i/dx_j)$ is not zero and so $q$ has a point where
it evaluates non-zero. On the other hand it is clear that $q$ vanishes
at any zero of $V$.

Conversely, suppose that $\{x_j\}$ are algebraically dependent. Then there is
some polynomial $q$ in $n$ variables with $q(\{x_j\})=0$. Consider
ideal $I = \langle f(\{d_i\},\{x_j\}),q(\{x_j\})\rangle$
and its variety 
$W$ in
$\mathbb{C}^{2n}$ (where we abuse notation with $d_1,\dots, d_n $ variables). 
This variety has dimension $n-1$ because it is isomorphic to
$V(\langle q(\{x_j\})\rangle ) $ in $\mathbb{C}^{n}$ under the isomorphism 
$(\{d_i\},\{x_j\}) \to \{x_j\}$. On the other hand if the elimination ideal
$ I  \cap \mathbb{C}[\{d_i\}]$
 is empty then 
it follows from the closure theorem (see Theorem 5.1) that
$W$
 has dimension at least $n$. This proves the existence of a non zero polynomial
$p(\{d_i\})$ in $I$. This polynomial  
evaluates to zero on the specific dimensions associated with the 
point $\{x_j\}$ since the generators of $I$ vanish on these
points.

Any algebraically independent set $\{x_j\}$ defines an algebraically independent
set $\{d_i\}$ for which $V$ is not empty. It follows that $V$ is not empty for all
algebraically independent $\{d_i\}$ because $V$ is empty only if the ideal of $V$
contains a constant element of the field $\mathbb{Q}(\{d_i\})$. Also, 
any zero of $V$ for
algebraically independent $\{d_i\}$  has algebraically independent coordinates
$\{x_j\}$ and so every point of $V$ has ~det$(J)$ non-zero. It follows that ~dim$(V)=0.$

\end{proof}
The next theorem is our variant of  
Laman's theorem.

\begin{thm}\label{T:Laman}
 Let $G$ be a maximally independent graph with 
$e$ edges and normalised
constraint equations $\{f_i\}$, and let $V$ 
be the associated variety in $\mathbb{C}^{e-1}$ for algebracialy independent $\{d\}$.
Then ~~dim$(V) = 0$.
\end{thm}

\begin{proof}
The normalised constraint equations have the form
required by Theorem 2.3 above while Theorem 6.5 of [9] implies that $\det(J)$ is not
zero as a polynomial in $\{x_j\}$.
\end{proof}

We shall use elimination theory to study
the varieties arising from various ideals generated by the constraint equations.
In order to keep track of the nature of solutions (whether they are
radical or not) it will be important, as we have intimated in the introduction,
to identify generators of one variable elimination
ideals which are irreducible polynomials.  Theorem \ref{T:Vb}
below will be needed to achieve this.

\begin{defn}\label{D:prime}
Let $I$ be an ideal in the polynomial ring $k[x_1,...,x_m]$
over a field $k$ of characteristic zero. Then
$I$ is {\it prime} if whenever $fg$ is in $I$ then either
$f$ is in $I$ or $g$ is in $I$.
\end{defn}

\begin{prop}\label{P:prime}
If $I$ is a prime ideal in $k[x_1,...,x_m]$ and if $\{x_{i_1},...,x_{i_t}\}$ is a
subset of
$\{x_i\}$ then the elimination ideal
\[
I \cap k[x_{i_1},...,x_{i_t}]
\]
is also a prime ideal.
\end{prop}

We now make  a simple but important observation. The
constraint equations for a graph are a parametric set 
when viewed as equations in the vertex 
coordinate variables and the dimensions. 
Indeed they are parametric in the vertex coordinate variables.
From this it follows that various associated complex
algebraic varieties are irreducible. 
For a discussion of such irreducibility 
see \cite{clos}. Thus we have the following general theorem which in turn
gives the irreducibility of what we call the big variety $V_b$.

\medskip
\begin{thm}\label{T:parametirc} Let $x = \{x_1, \dots ,x_m\}, d = \{d_1, \dots ,
d_r\}$ be indeterminates defining the polynomial ring
$\mathbb{Q}[x,d]$. Let $f_i(x,d)$ be polynomials of the form
$h_i(x) - d_i,  1 \le i \le r$, and let $I$ be the ideal of
polynomials in $\mathbb{Q}[x,d]$ which vanish on the variety
determined by $\{f_i: 1 \le i \le r\}$. Then $I$ is a prime ideal.
\end{thm}

\begin{thm}\label{T:Vb}
Let $G$ be a maximally independent graph with $n$ vertices 
and let $\{f\}$ be the
normalised constraint equations for $G$ for
the dimension set $\{d_e\} = \{d_1,...,d_r\}$ (where $r =
|E(G)|-1$). Let  $V_b \subseteq \mathbb{C}^{2n-4+r}$ be the
complex affine variety determined by $\{f\}$ as polynomial
functions belonging to
\[
\Q[d_1,...,d_r, x_1,...,x_{n-2},y_1,...,y_{n-2}].
\]
Then  $V_b$ is irreducible.
\end{thm}

\section{Connectedness and quadratic solvability}\label{S:qs}

The most tractable  CAD graphs from the perspective of solvability
are those which can be reduced to a collection of triangle
graphs by successive disconnections at vertex pairs. In this
section we indicate the way in which these graphs are 
quadratically soluble. We also recall
various  notions of  connectivity for graphs.

\medskip

\begin{defn}\label{D:qs}
Let $G$ be a maximally independent graph and let $V$ be the variety
 defined by  the constraint equations with generic normalised dimensions $\{d_e\}$.

(i) $G$ is said to be (generically) quadratically soluble (or simply QS) if
 every coordinate of every point
of $V$ lies in 
an extension of the base field 
$\mathbb{Q}(\{d_e\})$ of degree $2^n$ for some $n$.

(ii) $G$ is said to be soluble by radicals (or RS), or, simply, 
soluble, if
every such coordinate lies in a radical extension of the base
field.
\end{defn}

One could equally well define what it means for a specific
dimensioned graph to be QS or RS. For example it would be of
interest to know if particular graphs with integral dimensions are
soluble. Such problems lead rapidly into arithmetical problems
associated with multi-variable diophantine analysis and, with the
exception of some considerations of integral  doublets, we shall
not address such non-generic issues.

The field $\mathbb{Q}(\{d_e\})$ is the field of fractions of
polynomials in the dimensions. An irreducible quadratic polynomial
over this base field determines a field extension of degree 2 and
so a sequence of $n$ irreducible quadratic polynomials, with
coefficients in the new fields, give rise to a final field
extension of degree  $2^n$. Moreover any field extension of this
degree arises in this way. It
follows that if a maximally independent planar graph $G$ is
constructed through a sequence of triangles joined at
 common edges then $G$ is  QS. However,
 as is evident from Figure 3, not all QS graphs are
triangulated in this way.

\begin{figure}
\centering
\includegraphics[width=6cm]{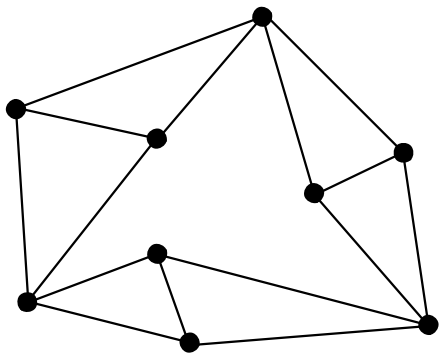}
\caption{A QS graph.}
\end{figure}

Recall that a graph $G$ is $n$-connected if there does not exist a
separation set with $n-1$ vertices. Thus the doublet is
$3$-connected while the graph of Figure 3 is $2$-connected. 
The following sufficient condition for quadratic solubility
was obtained in
Owen \cite{Owen}.

\begin{thm}\label{T:Owen} A CAD graph is (generically) QS if it admits a
reduction to triangle graphs by a process of
repeated separation at two-point
separation sets in which  all
but one of the separation components (the non rigid ones) have an edge
added between the separation pair.
\end{thm}

Note that the graph in Figure 3 can be reduced to a collection of triangles
in the manner of Theorem \ref{T:Owen}.
Graphs which are not
algorithmically reducible in this way of necessity possess 
a component which is 
$3$-connected. Thus the main theorem of the present paper provides a 
converse to Owen's theorem in the case 
of graphs with a planar embedding;  algorithmic reducibility of 
a planar CAD graph
is a necessary
condition to be (generically) QS or RS.
\section{3-connected maximally independent graphs}

We now embark on a graph-theoretic analysis of maximally independent,
3-connected, planar graphs. We shall prove the following main graph
reduction theorem.\\

\begin{thm} \label{T:4.1} Let $G$ be a $3$-connected, maximally independent,
planar graph with $|G|>6$. Then $G$ has either
\begin{itemize}\item[(i)] an edge which can be contracted to give a $3$-connected, maximally independent planar graph with $|G|-1$ vertices, or
\item[(ii)] a proper vertex-induced subgraph with three vertices of
  attachment which is maximally independent.
\end{itemize}
\end{thm}

We begin by stating some definitions and properties from graph theory.

The order of a graph $G$, denoted $|G|$ is the number of vertices in $G$. The
degree of a vertex $v$ in $G$, denoted $\deg(v)$, is the number of
edges of $G$ which are incident to $v$ or equivalently the number of
neighbours of $v$ in $G$. An edge joining vertices $x$ and $y$
is denoted by $(xy)$.

It is assumed throughout this section that all graphs $G$ have $|G|
\geq 2$  and if $H$ is described as a subgraph of $G$, then also $|H|
\geq 2$, unless it is explicitly stated otherwise. A vertex-induced
subgraph $H$ of $G$ has the additional property that if vertices $x$ and
$y$ are in $H$ and the edge $(xy)$ is in $G$, then the edge $(xy$) is also in $H$.

Let $H$ be a graph or a subgraph with $v$ vertices and $e$ edges. Define the
freedom number of $H$, written {\it free}$(H)$,
to be $2v-e-3$. A graph $G$
is independent if all its subgraphs $H$ have the property {\it free}$(H) \geq  0$. The graph $G$ is maximally independent if it is
  independent and {\it free}$(G) =  0$.

The graph $G\backslash e$ is the graph $G$ with the edge $e$
deleted. If $G$ is independent then $G\backslash e$ is also
independent and {\it free}$(G\backslash e) = $ {\it
free}$(G)+1$.

The graph $G/e$ is the graph obtained from $G$ by contracting the
edge $e$. This means that if the edge $e$ joins vertices $x$ and
$y$ then $G/e$ is obtained from $G$ by deleting the edge $e$,
merging the vertices $x$ and $y$ and reducing any resulting double
edges to single edges. Any such double edge must derive from a
$3$-cycle in $G$ that contains the contracted edge $e$. Thus
$|G/e| = |G| -1$ and if the edge $e$ is in a total of $c$ $3$-cycles
of $G$ then {\it free}$(G/e) = $ {\it free}$(G)+c-1$.

    An edge $e$ in an independent graph $G$ is said to be contractible if $G/e$
is independent and $free(G/e)=free(G)$. A necessary condition for $e$ to be
contractible is thus that it is in exactly one 3-cycle of $G$. However, this
condition is not sufficient as we show in Lemma 4.5 below.

If $H$ is a vertex-induced subgraph of $G$ then $G\backslash H$ is
the subgraph of $G$ induced by the vertices of $G$ that are not in
$H$. Here $|G\backslash H| < 2$ is not excluded. Thus $|G| =
|H|+|G\backslash H|$. The vertices of $H$ that have neighbours in
$G\backslash H$ are the vertices of attachment of $H$ in $G$. A
vertex-induced subgraph $H$ with $v$ vertices of attachment is
described as proper if $|H| > v$. An internal vertex of $H$ is a
vertex of $H$ that is not a vertex of attachment. An internal edge
of $H$ is an edge that joins to at least one internal vertex.

All vertices $v$ of a $3$-connected graph $G$ with $|G|>3$ have $\deg
(v) \geq 3$. The $3$-cycle is the only $3$-connected graph with
$|G|<4$. If $G$ is $3$-connected and $|G| > 3$ then any pair of
vertices in $G$ are joined by at least 3 paths which are
internally disjoint. We call such paths independent.

We shall say that a graph is planar if it has a planar embedding.
A planar embedding of a $2$-connected graph $G, |G| > 2$,  divides
the plane into disjoint regions called faces. One of these faces
includes the points at infinity. Each face is bounded by a cycle
of edges in $G$.

There are certain subgraphs whose occurrence is enough to ensure
that the graph resulting from an edge contraction 
is definitely not $3$-connected. The simplest of these consists of
a $3$-cycle connected into the remaining graph by exactly three
edges as shown in Figure 4. We call this subgraph the limpet.
\begin{figure}
\centering
\includegraphics[width=6cm]{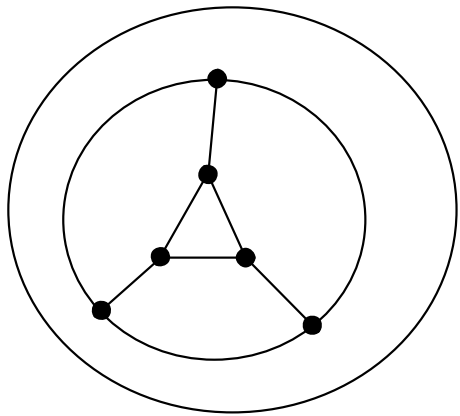}
\caption{The limpet subgraph}
\end{figure}
\phantom{xx}
 If a graph $G$ contains a limpet then $G$ also contains a
subgraph $H$ with three vertices of attachment in $G$, where $H$
is the subgraph induced by all vertices of $G$ that are not in the
$3$-cycle of the limpet. Clearly,  $|H| = |G|-3$ and $H$ has 6
less edges than $G$ so if $G$ is maximally independent then H is
also  maximally independent. If $|G| = 6$, then $G$ is the doublet.
If $|G| > 6$, then $H$ is a proper vertex-induced subgraph of $G$
with $3$ vertices of attachment that  is maximally independent.

The blocking role of the limpet should be clear by observing that
attaching the limpet by two vertices of attachment to any
contractible edge in a $3$-connected graph and assigning the third
vertex of attachment to any other vertex gives a $3$-connected
graph for which the result of contracting that same edge is
definitely not $3$-connected. This is shown in Figure 5. By adding
limpets into a graph in this way it is easy to generate graphs,
all of whose contractible edges produce graphs that are not
$3$-connected. Case (ii) of   Theorem \ref{T:4.1} is needed to
deal with limpets.

\begin{figure}
\centering
\includegraphics[width=8cm]{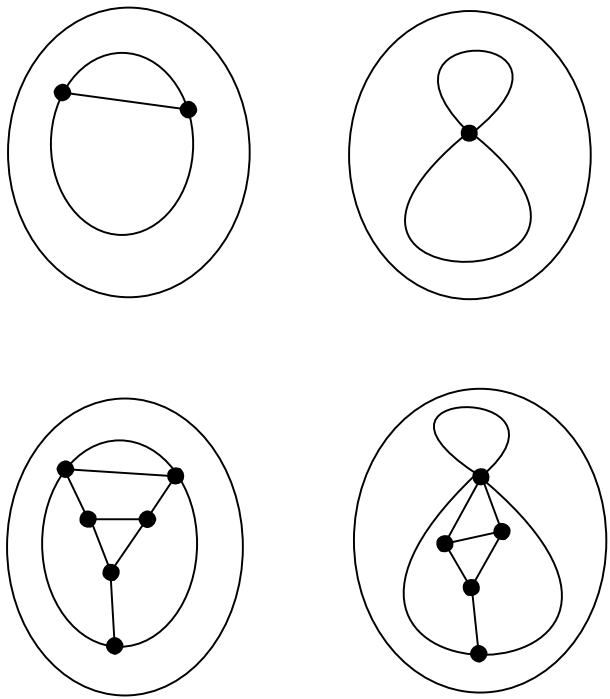}
\caption{Adding a limpet subgraph to an edge contraction.}
\end{figure}

We are now in a position to prove the main theorem of this section using the
sequence of lemmas proved below. To give some motivation to these lemmas we
begin with the proof of the main theorem.
\bigskip

{\it Proof of  Theorem \ref{T:4.1}.} Suppose that $G$ has no proper
vertex-induced subgraph with three vertices of attachment that is
maximally independent.

Assume for the sake of a proof by contradiction that $G$ contains no edge $e$
such that $G/e$ is $3$-connected and maximally independent.

$G$ is not the doublet because $|G|>6$ and $G$ has no limpets because it is
maximally independent and has no proper vertex-induced subgraph which is
maximally independent with three vertices of attachment.

 By Lemma \ref{L:4.7}, $G$ has no degree $3$-vertex on a $3$-cycle.

By the Corollary \ref{C:4.12}, $G$ contains an edge $e$ joining
vertices $x$ and $y$ such that $G/e$ is maximally independent.
Then $G/e$ is not $3$-connected, by the assumption, and by  Lemma
\ref{L:4.17}, $G$ has a $3$-vertex separation set $(x, y, w)$ for
some $w$, and this set separates $G$ into 2 proper components
$H_1$ and $H_2$. Let $H = H_1$ if $|H_1|<|H_2|$ otherwise $H =
H_2$. Now chose $e$ in $G$ which gives a minimal value for $|H|$.

By Lemmas \ref{L:4.17} and \ref{L:4.16} , the subgraph $H$
contains an edge $k$ which is internal to $H$ and which is
contractible as an edge in $G$. Thus $G/k$ is not $3$-connected by
the assumption. By Lemma \ref{L:4.18}, $k$ generates a $3$ vertex
separation set which has one proper component properly contained
in $H$. This contradicts the minimal condition on $|H|$
and completes the proof. \hfill $\Box$
\medskip

This proof requires a number of lemmas which deal with the effect
of an edge contraction on both maximal independence and
$3$-connectivity. The apparent complexity of the proof, including
the lemmas, is a result of the need to find edge contractions
which maintain both of these properties simultaneously.

The first three lemmas give some useful properties of maximally independent
graphs and subgraphs.

\begin{lma} \label{L:4.2} Let $H_1$ and $H_2$ be maximally independent
subgraphs of an independent graph $G$ with $|H_1 \cap H_2|  \geq
2$. Then $H_1 \cup H_2$ and $H_1 \cap H_2$ are both maximally
independent.
\end{lma}

\begin{proof}  $H_1 \cup H_2$ and $H_1 \cap {H_2}$ are both subgraphs of $G$ so
they are both independent. Let $H_1$, $H_2$, $H_1 \cup  H_2$ and
$H_1 \cap {H_2}$ have $v_1$, $v_2$, $v_u$, $v_i$ and $e_1$, $e_2$,
$e_u$, $e_i$ vertices and edges respectively. We have
$$2v_1-e_1-3=0, 2v_2-e_2-3=0, v_u=v_1+v_2-v_i, e_u =
e_1+e_2-e_i.$$
Thus {\it free}$(H_1 \cup H_2) = 2v_u-e_u-3 =
3-2v_i+e_i = -$ {\it free}$(H_1 \cap H_2)$.

Since both $H_1 \cup H_2$ and $H_1 \cap {H_2}$ are independent
they both have freedom numbers greater than or equal to zero and
thus equal to zero.
\end{proof}

\begin{lma} \label{L:4.3} Let $G$ be a maximally independent graph. Then $G$ is
$2$-connected.
\end{lma}

\begin{proof}
Suppose to the contrary.
Then there exist vertex-induced
subgraphs $H_1$ and $H_2$ such that
$G = H_1 \cup H_2$
and
$|H_1 \cap H_2| = 1$.
Using the same notation as for Lemma 4.2 we have
$$free(G) = 2v_u - e_u -3 \ge 2(v_1+v_2-1) -e_1 -e_2 -3 = 1,$$ 
which contradicts the fact
that $G$
is maximally independent.
\end{proof}

\begin{lma} \label{L:4.4} Let $G$ be a maximally independent graph. Then for any
edge $e$ the contraction $G/e$ has at most one separation vertex.
\end{lma}

\begin{proof}
 Suppose the edge $e$ joins vertices $(x,y)$ in $G$ which become the vertex
$w$ in $G/e$. Then any separation vertex of $G/e$ which is different from $w$ is
also a separation vertex of $G$ contrary to Lemma 4.3.
\end{proof}

The next lemma gives a useful criterion for an edge to be contractible.

\begin{lma}\label{L:4.5}  Let $G$ be an independent graph. An edge $e
= (xy)$ of $G$ is contractible if and only if
\begin{itemize}
\item[(i)]     $e$ is on exactly one $3$-cycle $(x,y,z$) of $G$ , and
\item[(ii)]    there is no maximally independent subgraph $R$ of $G$,
  $|R| \geq  3$, such that
$x$ and $y$ are in $R$ and $z$ is not in $R$.
\end{itemize}
The condition (i) can be replaced with weaker condition $(i')$ $e$
is on one or more $3$-cycles of $G$.
\end{lma}

\begin{proof} By definition e is contractible if and only if {\it
  free}$(G/e) = free(G)$ 
and $G/e$ is independent. We show that the first of these conditions is equivalent to
(i) and the second equivalent to (ii).

If $e$ is on $c$ $3$-cycles then {\it free}$(G/e) = free(G) + c+1-2$, so
{\it
  free}$(G/e) = free(G)$ if and only if $c=1$.

Now suppose (i) is true and (ii) is false. Then there is a
maximally independent subgraph $R$ of $G$ such that $x$,$y$ are in
$R$ and $z$ is not in $R$. We have {\it free}$(R) = 0$ and $R$
contains $e$, but not $z$. Thus {\it free}$(R/e) = -1$ (because
$R$ contains no $3$-cycle containing $e$) so $G/e$ is not
independent.

Conversely, suppose $G/e$ is not independent. Then $G/e$ contains
a subgraph, say $R/e$, with {\it free}$(R/e) = -1$ (since
contracting an edge reduces {\it free}$(H)$ by at most $1$ for any
subgraph $H$ of $G$). $R/e$ must contain the edge $e$ (or $R/e$
would also be a subgraph of $G$) so $R/e$ does indeed derive from
a subgraph $R$ in $G$ following contraction of $e$. Thus $R$
contains vertices $x$ and $y$ and {\it free}$(R) = 0$. The vertex
$z$ cannot be in $R$ because this would give {\it free}$(R/e) =
0$.

Clearly (i) implies $(i')$. Also $(i')$ and (ii) imply (i) because
 if $e$ is on two or more $3$-cycles then one of these contains a
 vertex $w$ different from $z$ and the $3$-cycle $(w,x,y)$ gives a subgraph $R$ which violates (ii).
\end{proof}

The next lemma is standard graph theory \cite{die} and describes
what happens if the result of an edge contraction in a 3-connected
graph is not $3$-connected.

\begin{lma} \label{L:4.6} Let $G$ be a $3$-connected graph. For any edge $e$
joining vertices $x$ and $y$, either $G/e$ is $3$-connected or $G$
has a $3$ vertex separation set consisting of $x$, $y$ and another
vertex $w$ of $G$.
\end{lma}

\begin{proof} Let $v$ be the vertex in $G/e$ that results from contracting
e and identifying x and y in G. If $G/e$ is not 3-connected then
it contains a separation pair $(a,w)$ and $a=v$ because $G$ is
3-connected. Thus $(v,w$) separate $G/e$ for some $w$ and
$(x,y,w)$ separates $G$.
\end{proof}

The next lemma identifies a class of $3$-connected 
independent graphs that always have a
contractible edge whose contraction gives a $3$-connected graph. These are
graphs that contain a $3$-cycle with one or two vertices with degree $3$.
Eliminating these graphs is helpful because the remaining graphs with a
$3$-cycle either contain a limpet or have all vertices on the 3-cycle with
at least two additional neighbours.

\begin{lma}\label{L:4.7} 
Let $G$ be a 3-connected, independent graph with no contractible
edges whose contraction gives a 3-connected graph. Then any 3-cycle in $G$
either has all its vertices with degree-3 or none of its vertices with
degree-3.
\end{lma}

\begin{proof}
 Suppose that $G$ contains a $3$-cycle $(x,y,z)$ with $\deg(x)
 = 3$. Let the third neighbour of $x$ be $t$. We will show that $\deg(y)= \deg(z)=3$.

We claim that both $(xy$) and $(xz)$ are contractible.

Suppose that neither $(xy)$ nor $(xz)$ is contractible. By Lemma
\ref{L:4.5} there is a maximally independent subgraph $R_{xy}$
containing $(xy)$ and not containing $z$ with $|R_{xy}| \geq 3$
and a maximally independent subgraph $R_{xz}$ containing $(xz$)
and not containing $y$ with $|R_{xz}| \geq 3$. By Lemma
\ref{L:4.3} the vertex $x $ has at least two neighbours in $R_{xy}
$which must be $y$ and $t$ and at least two neighbours in $R_{xz}$
which must be $z$ and $t$. Thus $R_{xy}\cap {R_{xz}}$ contains the
vertices $x$ and $t$ so $R_{xy} \cup R_{xz}$ is maximally
independent by Lemma \ref{L:4.2}. Then the subgraph $R_{xy} \cup
R_{xz}+(yz)$ has freedom number $-1$ (since $(yz)$ is in neither
$R_{xy}$ nor $R_{xz}$) which contradicts the independence of $G$.

Now suppose that $(xy)$ is contractible and that $(xz)$ is not.
Then $G/(xy)$ is not $3$-connected so there exists a separation
set $(x,y,w)$ of $G$. Since $G$ is $3$-connected each separation
component contains a vertex connected to $x$, so there are just
two separation components $C_z$ containing $z$ and $C_t$
containing $t$ and $w$ is distinct from $t$ and $z$. This is shown
in Figure 6. Then all paths from $x$ to $z$ in $G$ include the
edge $(xz)$ or include the vertex $y$ or include both the vertices
$t$ and $w$. If $(xz)$ is not contractible there exists maximally
independent $R_{xz}$ which includes $x$ and $z$ but not $y$. But
then all paths from $x$ to $z$  in $R_{xz}/(xz)$ include both
$t$ and $w$, so $t$ and $w$ are two separation vertices for
$R_{xz}/(xz)$ which contradicts Lemma \ref{L:4.4}.

We can now suppose that both $(xy)$ and $(xz)$ are contractible
and neither $G/(xy)$ nor $G/(xz)$ is $3$-connected. Then $G$ has a
separation set $(x,y,w)$ with a component $C_t$ which contains the
vertex $t$ and not the vertex $z$. $G$ also has a separation set
$(x,z,w')$ with a component $C_t'$ which contains the vertex $t$
and not the vertex $y$.

Since $t$ and $y$ are in different components of the separation
set $(x,z,w')$ all paths from $t$ to $y$ contain either $x$, $z$
or $w'$. The vertex set $(x,y,w)$ also separates $G$ and one
component $C_t$ contains $t$ (and not $z$) so there is a path from
$t$ to $y$ which lies inside $C_t$. Neither $z$ nor $x$ is inside
$C_t$ so $w'$ is in $C_t$ and $w'$ separates $y$ from $t$ inside
$C_t$. Since $G$ is $3$-connected this implies that the vertex $y$
is connected by the single edge $(yw')$ to $w'$ in $C_t$.
Similarly $w$ is in $C_t'$ and the vertex $z$ is connected by the
single edge $(zw)$ to $w$ in $C_t'$. This is shown in Figure 7
and Figure 8.

\newpage

\begin{figure}
\vspace{1in}
\centering
\includegraphics[width=7cm]{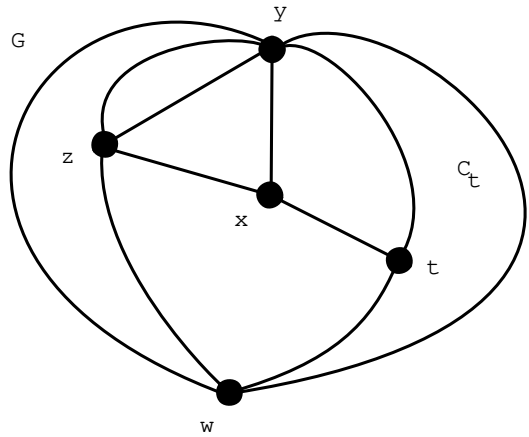}
\caption{$G/(xy)$ is  2-connected. The separation set $(x,y,w)$ in
$G$ gives two separation components.}
\end{figure}

\begin{figure}

\centering
\includegraphics[width=7cm]{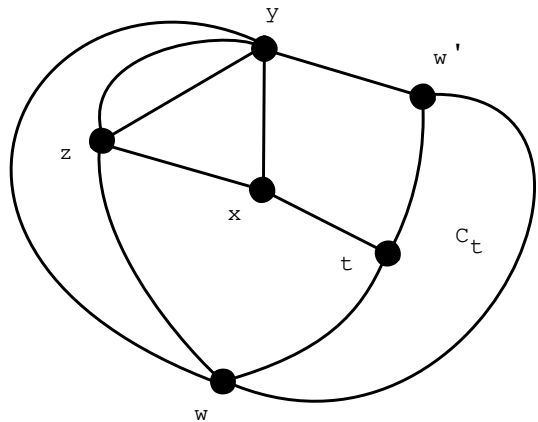}
\caption{If $(x,z,w')$ is also a separation set of $G$ then $w'$
is in $C_t$ and $w'$ is the only neighbour of $y$ in $C_t$}
\end{figure}

\phantom{xxx}
\newpage

\begin{figure}
\centering
\vspace{1in}
\includegraphics[width=10cm]{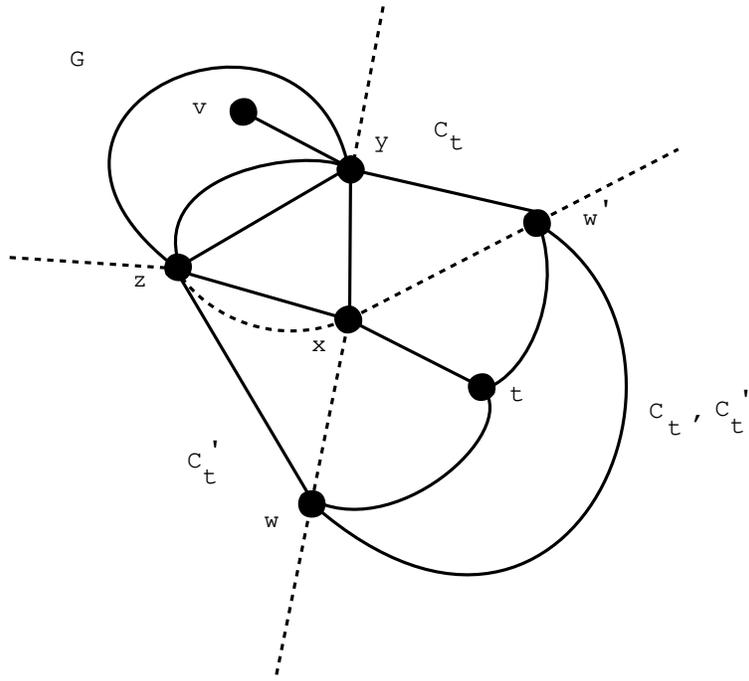}
\caption{Demonstration that $(y,z)$ is a separation pair of $G$ if
the contractions of both $(x,y)$ and $(x,z)$ are not 3-connected.}
\end{figure}

.

\pagebreak

Suppose that $y$ has a neighbour $v $ in addition to $x$, $z$ and
$w'$. Then $v$ is not in $C_t$ and $v$ is not in $C_t'$ because
$y$ is not in $C_t'$ and $v$ and $y$ are both distinct from
$(x,z,w')$. Since $t$ is in $C_t$ all paths from $v$ to $t$
include one on the separation set $(x,y,w)$ before any other
vertices of $C_t$. The vertex $x$ is connected only to $z$ outside
$C_t$ so a path including $x$ includes $z$. The vertex $w$ is
separated from $v$ by the separation set $(x,z,w')$ and of these
vertices only $z$ is outside $C_t$, so a path including $w$
includes $z$. Then all paths from $v$ to $t$ include either $y$ or
$z$, which contradicts the fact that $G$ is $3$-connected.

We conclude that $\deg(y)=3$ and similarly $\deg(z)=3$.
\end{proof}

The remaining lemmas make use of planarity in order to simplify
certain decompositions and to ensure a supply of contractible
edges. The first of these lemmas makes use of  the Kuratowski
theorem \cite{die} to simplify the number of separation components if the
result of contracting an edge is not $3$-connected.

\begin{lma} \label{L:4.8} Let $G$ be a $3$-connected, planar graph with a
$3$-vertex separation set. Then this separation set divides $G$
into exactly 2 proper components.
\end{lma}

\begin{proof} The separation set divides $G$ into at least 2 proper
components by definition. Suppose for a contradiction that there are 3
or more proper separation components. Then we can identify 3
vertices $w_1$, $w_2$, and $w_3$ each internal to a different
separation component. Let the separation set be the vertices
$v_1$, $v_2$ and $v_3$. There are paths connecting each of the
$w_i$ to each of the $vj$. By {Menger}$'$s theorem, the 3
paths from a $w_i$ to each of the three $v_j$ can be selected to
be internally disjoint because $G$ is $3$-connected and the paths
from different $w_i$ to any $v_j$ are internally disjoint because
they are in different separation components. Thus $G$ contains
$K(3,3$) as a topological minor contrary to { Kuratowski}$'$s
theorem.
\end{proof}

The next two lemmas lead to the Corollary \ref{C:4.12} that states
that every maximally independent, planar graph has at least 3
contractible edges. Lemma \ref{L:4.11} is stronger than is
required for this corollary but the greater detail will be useful
subsequently.

\begin{lma}\label{L:4.9}  Let $G, |G| > 2,$
be a $2$-connected planar graph with freedom number
$f$. Then every planar embedding of $G$ has the property
$$2(f-1)=\Sigma_i(n_i(i-4))$$
where the embedding has $n_i$ faces with
$i$ edges.
\end{lma}

\begin{proof} Let $G$ have $n$ vertices and $e$ edges and let the planar
embedding have $F$ faces. From { Euler}$'$s relation $F+n=e+2$ and
from the definition, $f=2n-e-3$ so $f=e-2F+1$. By definition $F=
{\Sigma_i(ni)}$. Each edge is in $2$ faces of the planar embedding
so $2e= \Sigma_i(i(ni))$ and the result follows by substituting
into $f= e-2F+1$.
\end{proof}

\begin{cor}
\label{C:4.10} A maximally independent, planar graph $G, |G| > 2$
contains at least one $3$-cycle.
\end{cor}

\begin{proof}  A maximally independent graph has $f=0$ and is $2$-connected by
Lemma \ref{L:4.3}. Thus in Lemma \ref{L:4.9} $n_3 \geq 2$, and the
boundary of one of these faces is a 3-cycle of $G$.
\end{proof}

\begin{lma} \label{L:4.11} Let $G$ be an
independent, planar graph which contains a
$3$-cycle $(x,y,z)$ and $(i,j,k)$ be any permutation
of $(x,y,z)$.

\begin{itemize}\item[(i)]     There exists a maximally independent
  subgraph $R_{ij}$ of $G$ with $i$ and $j$ in $R_{ij}$ and $k$ not in $R_{ij}$
  such that $R_{ij}$ contains an edge $e_{ij}$ which is contractible in $G$, and
\item[(ii)]    $R_{ij}\cap R_{jk} = j$.
\end{itemize}
\end{lma}

\begin{proof}
 Define the $R_{ij}$ as follows: If the edge $(ij)$ is
contractible then $R_{ij} = (ij)$. Otherwise, by Lemma \ref{L:4.5}
let $R_{ij}$ be a maximally independent subgraph containing $i$
and $j$ but not $k$ with $|R_{ij}|\geq 3$. Additionally take
$R_{ij}$ to be a maximal subgraph with these properties (maximal
in the sense that there is no subgraph $F$ with these properties
and $H \subseteq F$).

With this definition it is clear that $j$ is in $R_{ij}\cap
{R_{jk}}$. If $| R_{ij}\cap {R_{jk}}| \geq 2$ then {\it
}$free(R_{ij} \cup R_{jk}) = 0$ by {Lemma
  \ref{L:4.2}}. The vertices
  $i$ and $k$ are in  $R_{ij} \cup R_{jk}$ but the edge
$(ik)$ is not in $R_{ij} \cup R_{jk}$, so the subgraph $R_{ij}
\cup R_{jk} + (ik)$ of $G$ would  have freedom number  $-1$ which
contradicts the fact that $G$ is independent. Thus $| R_{ij}\cap
{R_{jk}}| = 1$ and $R_{ij}\cap {R_{jk}} = j$.

It remains to show that each $R_{ij}$ contains a contractible edge
which we do by induction. This is true for $|G|=3$. Assume it is
true for $|G| = N$.

Since every maximally independent planar graph contains a 3-cycle
(Corollary \ref{C:4.10}) it follows from the hypotheses that every
maximally independent, planar graph $R$ with $3 \le |R| \le N$ has
at least $3$ contractible edges. Thus if $(ij)$ is not
contractible then each $R_{ij}$ contains at least $3$ edges which
are contractible as edges in $R_{ij}$ and one of these, say edge
$e_{ij}$ is different from $(ij)$.

We claim that each $e_{ij}$ is also contractible as an edge in
$G$. Otherwise there exists a maximally independent subgraph $H$
in $G$, not contained in $R_{ij}$ but also containing  $e_{ij}$.
In fact $H\cap {R_{ij}} = e_{ij}$, because otherwise $H\cap
{R_{ij}}$ would be a maximally independent subgraph of $R_{ij}$
(by Lemma \ref{L:4.2}), containing $e_{ij}$ with $|H\cap {R_{ij}}|
\geq 3$ which contradicts the contractibility of $e_{ij}$ in
$R_{ij}$. Now $H \cup R_{ij}$ is also maximally independent by
Lemma \ref{L:4.2} and $|H \cup R_{ij}|
> |R_{ij}|$ which contradicts the maximality of $R_{ij}$ unless
$k$ is in $H \cup R_{ij}$. Suppose $k$ is in $|H \cup R_{ij}|$.
Then the independence of the subgraph  $H \cup R_{ij} + (ik) +
(jk)$ in $G$ requires $(ik)$ and $(jk$) in $H$ (since $k$ is  not
in $R_{ij}$). But $i$ and $j$ are in $R_{ij}$ and $H\cap  {R_{ij}}
= e_{ij}$ which would require $e_{ij} = (ij)$ contrary to the
assumption that $e_{ij}$ and the edge $(ij)$ are distinct.
\end{proof}

\begin{cor}
\label{C:4.12} Every maximally independent, planar graph $G$ has
at least 3 contractible edges.
\end{cor}

\begin{proof}
This was proved in Lemma \ref{L:4.11}.
\end{proof}

The next lemma guarantees the existence of a contractible edge in certain
subgraphs of an independent, planar graph.

\begin{lma}
\label{L:4.13} Let $H$ be a subgraph with 3 vertices of attachment
in an independent, planar graph $G$. If $H$ contains a $3$-cycle
with at least one vertex internal to $H$ then $H$ has an internal
edge that is contractible as an edge of $G$.
\end{lma}

\begin{proof} Let the $3$-cycle be $(x,y,z)$ with internal vertex $x$. By
Lemma \ref{L:4.11} there exist maximally independent subgraphs
$R_{xy}$ and $R_{xz}$ containing $(xy)$ and $(xz)$ respectively
and each of these contains a contractible edge.

We claim that either $R_{xy}$ or $R_{xz}$ have all their edges
internal to $H$. Otherwise both $R_{xy}$ and $R_{xz}$ each contain
at least two vertices of attachment, since if say $R_{xy}$
contains no vertex of attachment it is internal to $H$, and if it
contains one vertex of attachment then either all its edges are
internal to $H$ or $R_{xy}$ contains a vertex of $G\backslash H$. Then the
vertex of attachment would be a separating vertex for $R_{xy}$,
which contradicts Lemma \ref{L:4.3}. But if $R_{xy}$ and $R_{xz}$
each contain at least two out of the three vertices of attachment
then one of these vertices must be in both $R_{xy}$ and $R_{xz}$
and thus equal to the vertex $x$ since $R_{xy}\cap {R_{xz}} = x$.
This contradicts the requirement that $x$ is internal to $H$.
\end{proof}

The next sequence of lemmas has implications for $3$-connected
maximally independent planar graphs for which the contraction of
any contractible edge gives a graph which is not $3$-connected. We
have already shown that such a graph has a 3-vertex separation set
with exactly two components. The critical case for the proof of
theorem 4.1 is when each component has freedom number 1. The
difficulty is to show that each of these components contains a
$3$-cycle so that a reduction argument can be applied to the
smaller of the two components. Lemma \ref{L:4.9} alone is not
sufficient because substituting $f=1$ into this lemma leaves the
possibility that all faces have exactly $4$ edges. We exclude this
possibility by showing that at least one face has at least 5
edges.

\begin{lma}
\label{L:4.14} Let $G$ be a $3$-connected graph and let $H$ be a
proper vertex-induced subgraph of $G$ with 3 vertices of
attachment. If each vertex of attachment has at least $2$
neighbours in $H$ then $H$ is $2$-connected.
\end{lma}

\begin{proof}
 Suppose to the contrary that $H$ has a separation vertex $w$.
All three vertices of attachment cannot be in the same separation
component of $w$ because $G$ is $2$-connected. Thus there is a
separation component for $w$ which contains exactly one vertex of
attachment, say $v_1$ and this component must be just the edge
$(wv_1)$ or else $(w,v_1)$ would be a separation pair for $G$. This
contradicts the requirement that $v_1$ has at least $2$ neighbours
in $H$.
\end{proof}

\begin{lma}
\label{L:4.15}  Let $G$ be a $3$-connected planar, graph and let
$H$ be a proper vertex-induced subgraph of $G$ with $3$ vertices
of attachment and let each vertex of attachment have at least  $2$
neighbours in $H$. Then a planar embedding of $G$ implies a planar
embedding of $H$ and this embedding of $H$ has the three vertices
of attachment in one face boundary.
\end{lma}

\begin{proof}
$G$ has a planar embedding and deleting $G/H$ plus any edges
connected to $G/H$ gives a planar embedding of $H$. By Lemma
\ref{L:4.14} $H$ is $2$-connected, so the planar embedding of $H$
divides the plane into disjoint faces.

The three vertices of attachment of $H$ in $G$ are a separation
set for $G$. We claim that all vertices of $G/H$ lie in the same
face with respect to the embedding of $H$. By Lemma \ref{L:4.8}
the $3$-vertex separation set divides $G$ into exactly $2$
separation components. Thus every pair of vertices in $G/H$ is
joined together by a path in $G/H$.  All vertices of $G/H$ are
therefore embedded in the same face of the embedding of $H$
because otherwise these paths would cross a face boundary of the
embedding of $H$ and these face boundaries lie in $H$. There is a
vertex of $G/H$ adjacent to each of the three separation vertices
so the three separation vertices lie on this face boundary.
\end{proof}

\begin{lma}
\label{L:4.16} Let $G$ be a $3$-connected, independent, planar
graph and let $H$ be a proper vertex-induced subgraph of $G$ with
$3$ vertices of attachment $(v_1, v_2 \mbox{ and }  v_3)$ and let
each vertex of attachment have at least $2$ neighbours in $H$. If
$H$ has freedom number $1$ and if $H$ contains at most one of the
 edges $(v_1v_2)$, $(v_2v_3)$ or $(v_3v_1)$ then $H$ contains an edge
 adjacent to an interior vertex of $H$ that is contractible as an edge of $G$.
\end{lma}

\begin{proof}
 A planar embedding of $G$ gives a planar embedding of $H$. By
Lemma \ref{L:4.14} $H$ is $2$-connected and by Lemma \ref{L:4.15}
one of the face boundaries contains $v_1$, $v_2$ and $v_3$. Since
$H$ contains at most one of the edges $(v_1v_2)$, $(v_2v_3)$ or
$(v_3v_1)$ this face boundary has at least $5$ edges so by  Lemma
\ref{L:4.9} with $f=1$
 the embedding of $H$ has at least one face with $3$ edges and so $H$
 contains a $3$-cycle. Since $H$ has at most one of the edges
 $(v_1v_2)$, $(v_2v_3)$ or $(v_3v_1)$, $H$ has a $3$-cycle with an interior
 vertex and by Lemma \ref{L:4.13} $H$ contains an edge
adjacent to an interior vertex of $H$ that is contractible as an edge of $G$.
\end{proof}

\begin{lma}
\label{L:4.17} Let $G$ be a $3$-connected, maximally independent
 planar graph that contains no maximally independent vertex-induced
 subgraph with $3$ vertices of attachment and which has no degree $3$
 vertex on a $3$-cycle.
 For any contractible edge $e$
joining vertices $x$ and $y$, either $G/e$ is $3$-connected or $G$ has
 a $3$ vertex separation set consisting of $x$, $y$ and another vertex
 $w$ of $G$ with the following properties:
\medskip

1.      $G$ does not contain edges $(xw)$ or $(yw).$\\

2.      the separation set divides $G$ into exactly $2$ proper
components such that each proper component plus the edge $(xy)$ has
freedom number $1$.\\

3.      w has at least $2$ neighbours in each of the two proper components.\\
\end{lma}

\begin{proof}
Suppose $G/e$ is not $3$-connected. By Lemma \ref{L:4.6} and
\ref{L:4.8} $G$ has a $3$-vertex separation set $(x,y,w)$ which
separates $G$ into exactly $2$ proper components $C_1$ and $C_2$.
Let $H_1=C_1+(xy)+ (xw)'+(yw)'$ and $H_2=C_2+(xy)+ (xw)'+(yw)'$ ,
where $(xw)' = (xw)$ only if the edge $(xw)$ is in $G$ and
similarly for $(yw)'$. Let $G$, $H_1$ and $H_2$ have $v$, $v_1$,
$v_2$ and $e$, $e_1$, $e_2$ edges and vertices respectively. Let
$d=0$, $1$ or $2$ if none, one or both of $(xw)$ and $(yw)$ is in
$G$ and let $H_1$ and $H_2$ have freedom numbers $f_1$ and $f_2$.
We have
$$v=v_1+v_2-3, e=e_1+e_2-1-d, 2v-e-3=0,$$
$$
f_1=2v_1-e_1-3,  f_2=2v_2-e_2-3.$$
Thus $2(v_1+v_2-3)-(e_1+e_2-1-d)-3 = 0$ and so $f_1+f_2 = 2-d$

By hypothesis neither $H_1$ nor $H_2$ is maximally independent so
$f_1>0$ and $f_2>0$.  This requires $f_1=1$, $f_2=1$ and $d=0$.

Suppose a component, say $C_1$ has only vertex $a$ adjacent to
$w$. Then $H_1-w-(a w)$ has freedom number $0$ and $3$ vertices of
attachment in $G$. $H_1-w-(aw)$ is not the $3$-cycle because $w$
would be a degree $3$ vertex on a $3$-cycle contrary to hypothesis
so $H_1-w-(aw)$ is a proper maximally independent vertex-induced
subgraph of $G$, contrary to hypothesis.
\end{proof}

The final lemma allows us to conclude that under certain conditions one of
the separation components that can result from contracting an edge in a
subgraph must lie entirely within that subgraph.

\begin{lma}\label{L:4.18} Let $G$ be a $3$-connected graph and
let $H$ be a proper vertex-induced subgraph of $G$ with $3$
vertices of attachment $v_1, v_2 $ and $v_3 $ such that $G$ has
the edge $(v_1v_2)$ and does not have the edge $(v_2v_3)$ or the
edge $(v_1v_3)$ and let $v_3$ have at least $2$ neighbours in the
subgraph induced by the vertices of $G\backslash H+v_1+v_2+v_3$.
Then for any interior edge $e$ of $H$ either $G/e$ is $3$
connected or one of the separation components of $G/e$ is properly
contained in $H$.
\end{lma}

\begin{proof}
Let the edge e join vertices $x$ and $y$ with vertex $x$ interior
to $H$. Suppose $G\backslash e$ is not $3$-connected. By Lemma
\ref{L:4.6} $G$ has a $3$-vertex separation set $(x,y,w)$. See
Figure 9.

We claim that $w$ is in $H$. Suppose to the contrary that $w$ is
in $G\backslash H$. Since  $v_1$ and $v_2$ are adjacent they are internal
vertices of only one component so there is another component $C$
that has either none of $v_1$, $v_2$ or $v_3$ as an internal
vertex or contains $v_3$ and not $v_1$ and $v_2$ as internal
vertex. If $C$ contains none of $v_1$, $v_2$ or $v_3$ then there
is a path in $C$ from $w$ in $G\backslash H$ to $x$ in $H$ that
avoids all vertices of attachment contrary to the definition of
vertices of attachment. Suppose $C$ contains $v_3$ as an internal
vertex and not $v_1$ or $v_2$. The vertex $v_3$ has at least $2$
neighbours in $G\backslash H$ (because it has at least $2$
neighbours in $G\backslash H+v_1+v_2+v_3$ and $G$ does not contain
$(v_1v_3)$ or $(v_2v_3)$) so there is a vertex $u$ in $G\backslash
H$ that is a neighbour of $v_3$ and is different from $w$. See
Figure 10. Thus $u$ is in $(G\backslash H)\cap {C}$ and is
different from $v_1$, $v_2$, $x$, $y$ and $w$. One of the vertices
$v_1$ or $v_2$, say $v_1$ is not $x$ or $y$ and is thus in
$G\backslash C$. Now all paths from $u$ to $v_1$ include one of
$w$,$x$ or $y$ before any vertices in $G\backslash C$. All paths
in $C$ from $u$ to $x$ or $y$ contain $v_3$ and thus all paths
from $u$ to $v_1$ contain $v_3$ or $w$ contradicting the fact that
$G$ is $3$-connected.

Now $x$, $y$ and $w$ are in $H$ and one vertex of attachment, say
$v_1$ is different from $x, y$ and $w$. All vertices in $G\backslash
H$ are connected on paths excluding $x$, $y$ and $w$ so one
separation component contains at least $G\backslash H+v_1$ as
internal vertices and so the other component is properly contained
in $H$.
\end{proof}

\newpage

\begin{figure}
\vspace{1in}
\centering
\includegraphics[width=8cm]{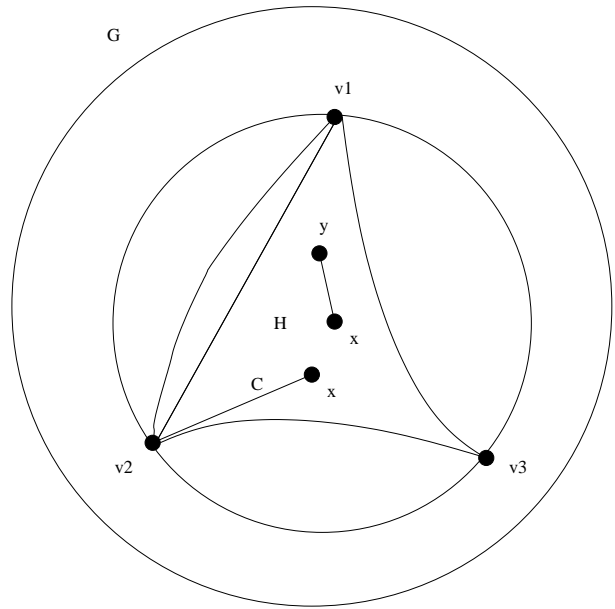}
\caption{The subgraph $H$ with 3 vertices of connection in $G$.
There are two different placings for an interior edge $e = (x,y)$ with
$x$ interior.}
\end{figure}

.

\newpage

\begin{figure}
\vspace{1in}
\centering
\includegraphics[width=8cm]{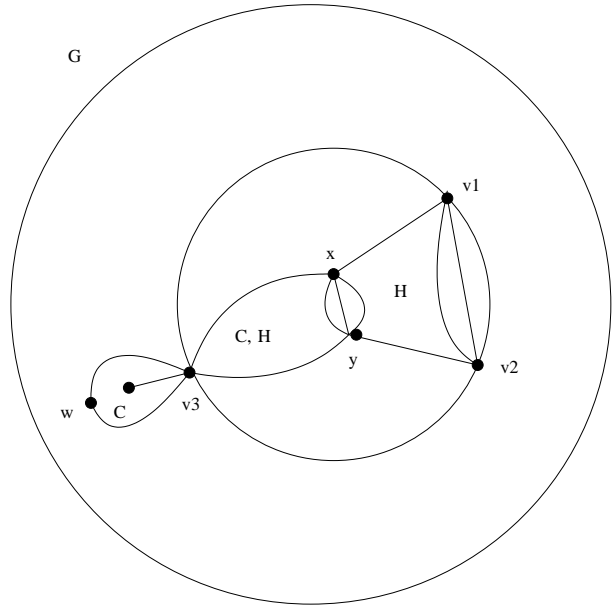}
\caption{The hypothetical structure of the separation
component $C$ if $\omega$ is in $G\backslash H$.
The vertex $y$ may be identical to $v_1$ or $v_2$.}
\end{figure}
.

\newpage

\section{Elimination ideals and
specialisation}\label{S:elimination}

In the present section we obtain irreducibility and
divisibility properties for generators of 
univariate elimination ideals and
their specialisations.
These properties play a prominent role in
the heart of our proof of the reduction step in that 
they connect the radical solvability
of generic equations with the radical solubility  
of the specialised equations.

Let $f_1,\dots,f_r$ be polynomials in the complex variables
$\{x_1,\dots,x_n\}$  which determine the complex algebraic
variety $V = V(f_1,\dots,f_r)$ in $\mathbb{C}^n$. For $ 1 \le t <
n$ the elimination ideal
\[
I_t = \langle f_1,\dots,f_r\rangle \cap \mathbb{C}[x_1,\dots,x_t]
\]
determines a variety $V(I_t)$ in $\mathbb{C}^t$. Plainly $V(I_t)$
contains $\pi_t(V)$, the projection of $V$ onto the subspace
$\mathbb{C}^t$. The following fundamental closure theorem may be
found in \cite{clos}.

\begin{thm}\label{T:closure}
The variety $V(I_t)$ is the Zariski closure of $\pi_t(V)$, that is, the
smallest affine variety containing $\pi_t(V)$.
\end{thm}

Let $\{d\} = \{d_1,\dots ,d_r\}$ be complex numbers forming an
algebricaly independent set with field extension
$\mathbb{Q}(\{d\})$.

\begin{thm}\label{T:generator}
 Let $\{ f \}$ be a set of polynomials in ${\mathbb Q} [d_1 ,
\ldots, d_r][ \{ x \}]$ which generates an ideal $I $
in $\mathbb{C}[\{x\}]$ 
 whose
complex variety $V(I)$ has dimension zero. Then each elimination
ideal \[ I_{x_i} = I \cap \mathbb{C}[x_i],
\] for
 $i = 1, \ldots, n$, is generated by a polynomial $g_i$ with
coefficients in ${\mathbb Q} [d_1 , \ldots, d_r]$ and $\deg(g_i) >
0$.  If, in addition, the set $\{ f \} $ generates a prime ideal
in the polynomial ring ${\mathbb
  Q}[d_1 , \ldots, d_r , x_1 , \ldots, x_n]$ then each $g_i$ may be
chosen to be irreducible in ${\mathbb Q} [d_1 , \ldots, d_r , x_i]$.
\end{thm}

\begin{proof}
Let $\hat{I}$ denote the ideal in ${\mathbb Q}(\{d_i\})[x_1,\dots ,x_n]$
generated by $\{f\}$ with elimination ideals
\[ (\hat{I})_{x_i} = \hat{I} \cap \mathbb{Q}(\{d\})[x_i].\]
Plainly, with the given inclusion  $ \mathbb{Q}(\{d\}) \subseteq \mathbb{C}$
we have  $\hat{I} \subseteq I$ and $I$ is the ideal in $\mathbb{C}[x_1,\dots ,x_n]$
generated by  $\hat{I}$. 

Since $(\hat{I})_{x_i}$ is an ideal in ${\mathbb Q} (\{d\})[x_i]$ it is
generated by a single polynomial $ g_i$, which is unique up to a
nonzero multiplier in ${\mathbb Q} (\{d\})$. Since $V(I)$ is nonempty
$g_i$ is not a nonzero constant, and so if deg$(g_i) = 0$ then
$g_i = 0$, and $(\hat{I})_{x_i} = \{0\}$. 
However, in this case we deduce that $I_{x_i} = \{0\}.$ This follows, for example,
from the fact that a basis for  $I_{x_i}$ may be derived from the generators
of $I$
by algebraic operations and so lie in $(\hat{I})_{x_i}$. (Consider a Groebner basis
construction for example.)
It now  follows that $V(I_{x_i}) = {\mathbb
C}$ and the closure theorem implies that the projection
$\pi_{x_i}(V)$ of $V(I)$ onto ${\mathbb C}_{x_i}$   is infinite and
hence that $V(I)$ is infinite, contrary to hypothesis. Thus deg$(g_i)
> 0$. 

The coefficients of $g_i$ are in ${\mathbb Q} (\{d\})$ and
so are ratios of polynomials in ${\mathbb Q}[\{d\}]$. Thus
we may replace $g_i$ by $p(d_1,...,d_r)g_i$ for some polynomial
$p$ to obtain the desired generator with polynomial coefficients. 
We may also arrange that the
highest common factor of the coefficients of $g_i$ is $1$.

We claim that the generator $g_i$, when viewed as an element of
the ring
${\mathbb Q} [\{d\},\{x_i\}]$, is also a generator for the
polynomial ring elimination ideal
\[
J_{x_i} = J  \cap {\mathbb Q} [\{d\},x_i],
\]
where $J$ is the ideal in $\mathbb{Q}[\{d\},\{x\}]$ generated by $\{f\}$.

Let $h \in J_{x_i}$. Then
$h$ is also in $(\hat{I})_{x_i}$ and so $h = qg_i$ with $q$ in  ${\mathbb
Q} (\{d\})[x_i]$. Clearing the denominators of the coefficients of
$q$ obtain the factorisation $r(d_1,...,d_r)h =
(r(d_1,...,d_r)q)g_i$ where $r(d_1,...,d_r)$ is in ${\mathbb Q}
[\{d\}]$  and $r(d_1,...,d_r)q$ is in ${\mathbb Q} [\{d\}][x_i]$.
Since, by the hypotheses, the  ideal
$J$ is prime, so too is $J_{x_i}$ and so one of these factors
belongs to $J_{x_i}$. However, if $r(d_1,...,d_r)q$ belongs to
$J_{x_i}$ then we can repeat the factorisation argument with
$r(d_1,...,d_r)q$ in place of $h$. Factoring in this way 
at most finitely many times
we see that we can 
assume that $h$ has the form  $pg_i$ with $p$ in ${\mathbb Q} (\{d\})[x_i]$.
Since the coefficients of $g_i$ have no common factor it follows
that $p$ is in ${\mathbb Q} [\{d\}][x_i]$ and that $g_i$ is a
generator for $J_{x_i}$. Since  
$J_{x_i}$ is prime  this in turn entails that the
generator $g_i$ is irreducible in ${\mathbb Q} [\{d\},\{x_i\}]$.
\end{proof}

We now show that in the case $r = 0$ the  specialised generator $g(d',x_i)$
is non-zero
and divisible by the generator $g_i(x_i)$ of the elimination ideal
of the specialised ideal. As we note below,
such divisibility may fail for a double specialisation !

For later convenience the role of $\mathbb{Q}$ in the theorem above is played 
below by $\mathbb{E} \subseteq \mathbb{C},$ a finite transcendental field extension of
$\mathbb{Q}$. (It is trivial to generalise the theorem above with
$\mathbb{Q}$ replaced by $\mathbb{E}$.) Specialisation occurs for the single
variable $d$ associated  with the transcendental extension $\mathbb{E}(d)$.
For an ideal  $\hat{I}$ in $\mathbb{E}[d][x_1,...,x_n]$ we shall write
$(\hat{I})'$ for the specialisation of $\hat{I}$ resulting from the substitution
$ d \to d'$.

\begin{thm}\label{T:AG4a}
 Let $\{ f \}$ be a set of polynomials in ${\mathbb E}
[d][x_1 , \ldots, x_n]$  which generate an ideal $\hat{I}$ 
in ${\mathbb E} (d)[x_1 ,\ldots, x_n]$ and an ideal $I$ in
$\mathbb{C}[x_1 , \ldots, x_n]$
 whose complex
variety $V(I)$ has dimension zero. Let $d' \in {\mathbb Q}$ be a
specialisation of $d$ giving rise to the set $\{ f' \}$ in
${\mathbb E}[x_1 , \ldots, x_n]$ with ideal 
$I'$ whose complex variety also has dimension zero.

Let $g(d, x_1)$  in ${\mathbb E} [d][x_1]$
 and $g'(x_1)$ in ${\mathbb
  E} [x_1]$ be generators for the elimination ideals
  $I_{x_1}$ and $(I')_{x_1}$
respectively, 
as provided by the previous theorem.
Finally, assume that the ideal in $\mathbb{E}[d,x_1,\dots,x_n]$
generated by $\{f\}$ is prime.
  Then

  (i) the specialisation $((\hat{I})_{x_1})'$
of $(\hat{I})_{x_1}$, is contained in $((\hat{I})')_{x_1}$,

(ii)  the
  degree of $g(d', x_1)$ is greater than zero, and

  (iii) $g'(x_1)$ divides
  $g(d', x_1)$.
\end{thm}

\begin{proof}
We have
\[
((\hat{I})_{x_1})' = \{p(d',x_1) : p \in \hat{I} \cap {\mathbb E}
(d)[x_1]\}.
\]
But if $p(d,x_1) \in \hat{I}$ then $p(d',x_1) \in (\hat{I})'$
and so $((\hat{I})_{x_1})' \subseteq ((\hat{I})')_{x_1}$.
Thus if $g(d',x_1)$ is not the zero polynomial then $g'(x_1)$ divides
$g(d',x_1)$ and deg($g(d',x_1)) > 0$.

Let $J$ be the ideal in ${\mathbb E} [d,x_1,...,x_n]$ generated by
$\{f\}$ and let $J_{d,x_1}$ be the elimination ideal $J\cap
{\mathbb E}[d,x_1]$. Then $J_{d,x_1}$ has generator $g_1(d,x_1)$
where this polynomial is the generator of $(\hat{I})_{x_1}$ in ${\mathbb
E} (d)[x_1]$ provided by the previous theorem.
By this theorem  we may assume that $g_1(d,x_1)$ is
irreducible in ${\mathbb E} [d,x_1]$. In this case it is not
possible to have $g_1(d',x_1) = 0$ for all $x_1$, for otherwise
$g_1$ would have a proper factor $(d-d')$.
\end{proof}

It is instructive to note that Theorem 5.3 is not valid without
the assumption that the big ideal is prime. Consider the equation
set
$$
(dx-1)p(x) = 0, ~~~~d(dx  - 1) = 0,
$$
where $p(x)$ is a polynomial in one variable $x$ over $\mathbb{Q}$
and $d$ is a single parameter. For generic $d$ the ideal $I =
\langle (dx-1)p(x),d(dx-1)\rangle$ in $\mathbb{Q}[x]$ is the
principal ideal $\langle dx-1\rangle, V(I)$ is the singleton
$\{1/d\}$ and ~dim$(V(I)) = 0$. For the specialisation $d=0$ the
ideal for the specialised equations is $I'=\langle p(x)\rangle$
and $V(I')$ is the finite set of zeros of $p$ and so is also zero
dimensional. However, it is not possible to choose a generator for
$I'$ which divides a nonzero generator of $I$, and so the
conclusion of Theorem 5.3 cannot hold for this equation set.

Note also that in this example we may choose $p(x)$ to be a
polynomial which is not soluble over $\mathbb{Q}$ so that while
the generic variety $V(I)$ is radical the variety for the
specialised equations is not radical.

It is also instructive to note that Theorem 5.3 is not valid for
the specialisation of more than one parameter. For example, let
\[f_1 = x_1(1-x_1x_2) - d_1, f_2 = x_2(1-x_1x_2)-d_2,
f_3 = x_3(1-x_1x_2)-d_3.
\]
For the double specialisation $d_1=d_2=0, V(I')$ is the single
point $ x_1=0,x_2=0,x_3=d_3$ and
\[
g_1(d_1,d_2,x_1)=d_2x_1^3 - d_1x_1 +d_1^2
\]
which becomes zero on this specialisation.
\section{The reduction step}\label{S:reduction}

Equipped with the elimination theory of the last section we are
now able to prove the reduction step stated in the introduction.

Let $G$ be a maximally 
independent graph with $n$ vertices and $r+1$ edges and suppose that
$G$ has an edge contraction to a maximally independent graph $G/e$. We label the
vertices so that $e$ is the edge $(v_{n-1}v_{n-2}), e$ is in the 3-cycle
$(v_n,v_{n-1},v_{n-2})$ and we regard $(v_{n-1}v_n)$ as the base edge. Furthermore, we
normalise the constraint equations $\{f_e\}$ so that the coordinates for the base
vertices are $(x_{n-1},y_{n-1}) = (0,0), (x_n,y_n) = (1,0)$. Let us label edges  so
that the contractible edge e is the rth edge, with the associated (squared)
dimension $d_r$, and the edge $(v_{n-2}v_n)$ has dimension $d_{r-1}$. Finally let
$f_1,...,f_r$ be a listing of the normalised constraint equations for $G$
compatible with this notation.

Now consider a set of normalised constraint equations for the
contracted graph $G/ e$. We lose two edges from $G$ (edge $r-1$
and edge  $r$) and we can take the normalised constraint equations
to be the equations $f_1,\dots f_{r-2}$ with the substitution
$x_{n-2} = 0$, $y_{n-2} = 0$.

First  consider the dimensions $\{d\} = \{d_1,\dots,d_{r-2}\}$
(together with $d_b = 1$) to be a generic set of real numbers.
Since the contracted graph is maximally independent the solutions
(for $x_1,\dots,x_{n-3},y_1,\dots,y_{n-3}$) form a zero
dimensional variety, $V(0,0)$ say. (The choice of notation will
become clear shortly.) Clearly this is essentially the variety of
the constraint equations $\{f_1,\dots,f_r\}$ for the
dimension set
\[
\{d_1,\dots,d_{r-1},d_{r-2},1 ,0\}
\]
for $G$
resulting from the double specialisation $d_{r-1} = 1, d_r = 0.$
Thus, in order to establish the reduction step it will be sufficient
to show that if $G$ is generically radical then
the variety arising from the semi-generic
double specialisation is also a radical variety.
This requires some care in view of the failure
of a double specialisation variant of the
Theorem \ref{T:AG4a}.
We shall break the double
specialisation into two steps. Also, instead of specialising
the generic edge lengths $d_r,d_{r-1}$ we choose to start afresh
and specialise the given coordinates $x_{n-2},y_{n-2}$. This
results in a simpler comparison of varieties.

In fact we can prove the reduction step for general non-planar
graphs.

\begin{thm}\label{T:reduction}
 Let $G$ be a  maximally independent  graph
which has an edge contraction to a  maximally
independent graph $G/ e$. If  $G$ is radically soluble then the
graph $G/ e$ is also radically soluble.
\end{thm}

\begin{proof}
Consider the set of dimensions $\{d\} = \{d_1,\dots,d_{r-2}\} $
and the constraint equations $\{f_1,\dots ,f_{r-2}\}$ in the
variables $x_1,\dots x_{n-3}$, $y_1,\dots,y_{n-3}$ which arise
when the pair $(x_{n-2},y_{n-2})$ takes three possible pairs of
values, namely $(X,Y), (X,0)$ and $(0,0)$, where $X,Y$ are
generic. Denote the three corresponding "big" varieties, where
$\{d\}$ is a set of variables, by $V_b(X,Y))$, $V_b(X,0)$ and
$V_b(0,0)$. For generic values of $\{d\}$ let the corresponding
"small" varieties be $V(X,Y))$, $V(X,0)$ and $V(0,0)$. Also we
write $I_b(X,Y), I(X,Y)$ etc., for the six corresponding ideals

We have the following:

1. The varieties $V_b(X,Y))$, $V_b(X,0)$ and $V_b(0,0)$ are
irreducible. This follows from the fact that the equations are
parametric in the variables. See Theorem \ref{T:Vb}.

2. 
The variety $V(0,0)$ is zero dimensional by Theorem 2.4 because it is the
variety of the maximally independent generic graph $G/e$. The varieties $V(X,0)$
and $V(X,Y)$ also have the form required for Theorem 2.3. The determinant of
the Jacobian matrix for $V(0,0)$ is obtained from the corresponding
determinants for $V(X,0)$ and $V(X,Y)$ by substituting $X=0$ and $Y=0$ and thus
neither of the determinants of the Jacobian matrices for $V(X,0)$ and $V(X,Y)$
are identically zero. Then $V(X,0)$ and $V(X,Y)$ are zero dimensional by Theorem
2.3.

We may now apply the specialisation theorem of Section
\ref{S:elimination} two times, once for the specialisation $(X,0) \to
(0,0)$ and once for the specialisation $(X,Y) \to (X,0)$.

Suppose then, that $V(0,0)$ is non-radical. In fact assume that
there is a point of this variety whose $x$-coordinate is not in a
radical extension of  $\mathbb{Q}(\{d\})$. Since $V_b(0,0)$ is
irreducible and $V(0,0)$ is zero dimensional, it follows from
Theorem \ref{T:generator} that there exists a univariate
polynomial $g(x_i)$ in $\mathbb{Q}(\{d\})[x_i]$ which generates
the elimination ideal $I(0,0)_{x_i}$. By  the closure theorem,
Theorem \ref{T:closure}, $\pi_{x_i}(V(0,0))$ is precisely the
variety of the elimination ideal for $x_i$ and this is precisely
the set of zeros of $g_i$. By the non-radical hypothesis there
exists an $x_i$ such that $g_i$ has some of its roots non-radical
(over $\mathbb{Q}(\{d\})$). By irreducibility, all the roots are
non-radical.

Likewise, $V(X,0)$ is zero dimensional and there exists a
polynomial $g(x_i,X)$, with positive degree in $x_i$, which
generates $I(X,0)_{x_i}$. Moreover, since $V_b(X,0)$ is
irreducible we may choose $g$ so that $g(x_i,X)$ is not divisible
by $X$  and hence $g(x_i,0)$ is not identically zero. But
$g(x_i,0)$ is in $I(0,0)_{{xi}}$ and  so $g(x_i) $ divides
$g(x_i,0)$. Thus $g(x_i,0)$ has a non-radical root, $g(x_i,0)$ is
non-radical and $V(X,0)$ is non-radical.

Repeating this argument for $V(X,0)$  and $V(X, Y)$ shows that
$V(X,Y)$ is non-radical over $\mathbb{Q}(\{d\})$. Thus $V$ is
non-radical over $\mathbb{Q}(\{d\},X,Y).$ However, by triangle
geometry $X$ and $Y$ are radical functions of $d_{r-1}$ and $d_r$.
Thus  $V $ is non-radical over $\mathbb{Q}(d_1,...,d_r).$
\end{proof}
\medskip

\noindent {\bf Remark.} One needs to take care with simultaneous
specialisation. If we do both specialisations together on $V(X,Y)$
we might have
\[
g(x_i, X,Y) = Xp(x_i,X ,Y)+Yq(x_i,X ,Y), \]
 where, for example, $Y$  does
not divide $p$ and so $g(x_i, 0, 0) = 0,$ which gives no
information on divisibility. In fact we have not excluded this
possibility by doing the specialisations one at a time. However we
have shown that if this does occur then $p$ and $q$ both have
factors which are non-radical. This is sufficient to  deduce that
$g(x_i, X, Y)$ is non-radical, even if it is zero on the double
specialisation.
\section{Galois group under specialisation}

We now obtain a theorem concerning the Galois groups of polynomials whose
coefficients contain indeterminates which may be specialised.
This theorem plays a role in the proof of the fact 
that if the graph $G$ is soluble by
radicals for generic dimensions then  it is also soluble by radicals
for certain specialised dimensions. 
In the proof we
 make use of the identification of the Galois group of $p$ as the
set of permutations in an 
index set associated with a certain irreducible factor of a
multi-variable polynomial constructed from $p$. This
identification is well-known and given in Stewart \cite{ste}.

Let ${d} = \{d_1,\dots,d_n\}$ be algebraically independent
variables with the rational field extension $\mathbb{Q}({d})$ and
let ${d'} = \{d_1',\dots,d_n'\}$ be an $n$-tuple of rationals,
viewed as a specialisation of ${d}$. 

\begin{thm}\label{T:galois1}
Let $p \in   {\mathbb Q}  [{d}]  [t]$ be an irreducible monic
polynomial with Galois group ~~Gal$(p)$ when viewed as a
polynomial in ${\mathbb Q}   ({d})  [t]$. Let ${d}' \in {\mathbb
Q}^n$ be a specialisation of ~${d}$ and let ~$p'$ be the
associated specialisation of $p$ with Galois group Gal$(p')$ over
${\mathbb Q}$. Then Gal$(p')$ is a subgroup of Gal$(p)$. In
particular if $p$ is a radical polynomial then so too is $p'$.
\end{thm}

\begin{proof}
 Consider the irreducible polynomial
 $$ p(t) = t^m + b_{m-1}({d}) t^{m-1} +
\ldots + b_0 ({d})$$ with coefficients  $b_i ({d})$ in $ {\mathbb
Q} [{d}]$.  Let $\alpha_1, \ldots, \alpha_m$ be the roots of
$p(t)$  in some splitting field, let $\{x_1, \ldots , x_m\}$ be
indeterminates  and let
$$\beta = \alpha_1 x_1 +
\ldots + \alpha_m x_m.$$
Let $S_m$ be the symmetric group  and define
$$
Q(t, x_1, \ldots , x_m) = \prod_{\sigma \in S_m} \left(t - \sigma
  (\beta) \right)
$$
where $\sigma (\beta) = \alpha_1 x_{\sigma (1)} + \ldots +
\alpha_m x_{\sigma (m)}$.  On expanding the product it can be seen
that the coefficient of a monomial $t^k{x_1}^{{i_1}}{x_2}^{{i_2}},
\dots ,{x_m}^{{i_m}}$ is a symmetric polynomial in the roots
$\alpha_i$. It follows that these coefficients are
 polynomials in $b_{m-1}
({d}), \ldots, b_0 ({d})$. (See \cite{ste}.) Thus the polynomial
$Q$ belongs to
${\mathbb Q} [{d}][t, {x}].$

Let $Q = Q_1 Q_2 \ldots Q_r$ where each $Q_i$ is irreducible in
${\mathbb Q} [{d}] [t, x]$ and where $Q_1$ contains the factor $(t
- \beta)$.   Since the roots of an irreducible polynomial are
distinct so too are the expressions $\sigma(\beta)$ and it follows
that the polynomial  $Q_1$ is well-defined.

We have
$$Q_1 = \prod_{\sigma \in S} (t - \sigma ( \beta ))$$
for some index set $S$. This index set is a subgroup of $S_m$
which is identifiable with the Galois group of $p$. It coincides
with the group of permutations $\sigma$ of the variables $x_1
\ldots, x_m$ for which $\sigma (Q_1) = Q_1$. In fact each $Q_i$
has the form $\tau(Q_1)$ for some permutation $\tau$ and from this
it follows that if $\sigma(Q_i) = Q_i$ for some $i$ then this
holds true for all $i$ and $\sigma$ is in the Galois group.

Now consider the specialisation $Q'$ of the polynomial $Q$ in $ ~
{\mathbb Q} [{d}] [t, {x}]$ upon replacing ${d }$ by $ {d'}$.
Since the coefficients of $Q$ are polynomials in $b_{m-1} ({d}),
\ldots, b_0 ({d})$ it is easy to see that $Q'$ coincides with the
'$Q$ polynomial' for $p'$. Thus $Q'$ is equal to the polynomial
$$\prod_{\sigma \in S_{m}} ~ \left(t- \sigma ( \beta') \right)$$
where $\beta' = \alpha_1' x_1 + \ldots + \alpha_m' x_m$ and
$\alpha_1' , \ldots, \alpha_m'$ are the roots of the
specialisation $p'$ in some order. (Despite the notation we do not
imply that there is a link between any $\alpha_i'$ and
$\alpha_i$.)

Note that for any permutation $\sigma$ and polynomial $P$ in ${\mathbb
  Q} [{d}][t, {x}]$ the polynomial $\sigma(P)$ is
  defined by permuting the indeterminates $x_1, \ldots, x_m$.  Thus
  $\sigma(P)' = \sigma (P')$, which is to say that
  the permutation action on these
  polynomials commutes  with specialisation.

Consider now both the specialisation of the factorisation, namely
$$Q' = Q_1' Q_2' \ldots Q_r',$$
and the irreducible factorisation of $Q'$ in ${\mathbb Q} [t,
{x}]$, namely
$$Q' = P_1 P_2 \ldots P_s.$$

Let us assume first that the roots $\alpha_i'$ are distinct. Then,
since each $P_i$ is necessarily a product of some of the
irreducible factors $t - \sigma (\beta')$, there is a unique
factor, $P_1$ say, divisible by $t - \beta'$.  Once again (and
even though $p'$ may be reducible) the Galois group $Gal (p')$ is
identifiable with $T$ where $T \subseteq S_m$ is the index set
such that
$$P_1 = \prod_{\sigma \in T} (t - \sigma (\beta')).$$
The roots $\alpha_1', \ldots, \alpha_m'$ do not correspond to
$\alpha_1, \ldots, \alpha_m$ and so we cannot assume that $P_1$
divides $Q_1'.$ (Such divisibility  gives $T \subseteq S$ and so
completes the proof in this case.) However, let $\sigma \in T$, so
that $\sigma (P_1) = P_1$, and suppose that $P_1$ divides $Q_i'$.
Then $P_1$ divides $\sigma (Q_i') = \sigma (Q_i)' = Q_j'$ say,
where $Q_j = \sigma(Q_i)$.  By the distinctness of the roots
$\alpha_i'$ and the fact that $\mathbb{Q}[t,x]$ is a unique factorisation
domain, it follows that if $P_1$ divides both $Q_i'$ and
$Q_j'$ then $i=j$. Thus $\sigma (Q_i) = Q_i$.  But by our remarks earlier
this condition
on $\sigma$ is equivalent to $\sigma(Q_1) = Q_1$ and hence $\sigma
\in S = Gal (p)$.

We now give more notational detail on this case which we shall
elaborate further to prove the general case.

Assume that $p' = h_1 h_2 \ldots h_q$ where $h_1, \ldots, h_q$ are
distinct irreducible polynomials in ${\mathbb Q} [t]$ with $\deg
h_i = r_i$.

 The Galois
group $T = Gal(p')$ can be identified in a natural way with a
subgroup of the product group $Gal(h_1) \times \dots \times
Gal(h_q)$. We remark that $T$ may be a proper subgroup. For
example, if $h_1$ and $h_2$ determine the same field extension of
$\mathbb{Q}$ then $r_1 = r_2$ and $Gal(h_1h_2) = Gal(h_1)$. (Each
permutation of the roots of $h_1$ determined by an element of
$Gal(h_1)$ is matched with a corresponding permutation of roots of
$h_2$.) In general $Gal(p')$ is a product of the Galois groups of
the distinct field extensions determined by irreducible factors of
$p'$.

The irreducible polynomial $P_1$ above factors as a product
$$P_1 = \prod_{
\sigma = \sigma_{1} \times \ldots \times \sigma_{q} \in T} \left(t
- (\sigma_{1} (\beta_1') + \ldots + \sigma_{q} (\beta_{q}')
\right),$$ where $\beta_i' = \alpha_{i,1}' x_{i,1} + \ldots +
\alpha_{i,
  r_{i}}' x_{i, r_{i}}$, and where $\alpha_{i,1}' , \ldots, \alpha_{i_{i}
  r_{i}}'$ are the distinct roots of $h_i$.
   Thus we have $r_1+ \dots + r_q = m$ and we have  identified the
variables $x_{1}, \ldots, x_m$ with the variables
$$x_{1, 1}, \ldots, x_{1,r_{1}}, x_{2,1},\dots  ,x_{2,r_2},~~~
  \ldots ~~~ , x_{q,1}, \dots ,x_{q, r_{q}}.$$

Consider now the general case wherein
$p' = h_{1}^{n_{1}} h_{2}^{n_{2}} \ldots
h_{q}^{n_{q}}$ where each $h_{i}$ is as before, with degree $r_i$.
Now each root $\alpha_{i,k}'$ appears with multiplicity $n_{i}$
and $m$ now satisfies the equation $$n_1r_1+\ldots + n_qr_q = m.$$
Let us accordingly relabel the variables $x_{i,j}$ as
$$\{ x_{i,k,t} : 1 \leq i \leq q, 1 \leq k
\leq r_i, 1 \leq t \leq n_i \}$$
Identify each element $\sigma = \sigma_1 \times \ldots \times
\sigma_q$ of $T = Gal(p')$ with the permutation in
$$ (Gal(h_1)\times \ldots \times Gal(h_1)) \times \ldots \times
(Gal(h_q)\times \ldots \times Gal(h_q))$$ which respects the
ordering of repeated roots and which respects the matching of
permutations in $Gal(h_i)$ and $Gal(h_j)$ if $h_i$ and $h_j$
determine the same field extension. In this way we obtain an
identification of $Gal(p')$ as a subgroup of $S_m$. Note that
there is a degree of choice in this identification; the
permutations that permute only indices of equal roots give rise to
distinct embeddings.

 Consider now
the  polynomial in $\mathbb{Q}[t,x]$
associated with this inclusion defined by
$$P_{*} = \prod_{\sigma \in Gal(p')\subseteq S_m} \left(t - \sigma (\beta') \right).$$
This polynomial  has the form $\hat{P}_1$ where
$$\hat{P}_1 = P_1(X_1,  \ldots  , X_q)$$
where $P_1$ is the irreducible polynomial we had in the previous
case and where each $X_i$ is the sum of
those variables corresponding to
repeated and matched roots.

Since $P_1$ is irreducible it follows
that $\hat{P}_1$  is irreducible.  It follows further that the
irreducible factors of $Q'$, and hence $Q_1'$, have the form
$\tau(\hat{P}_1)$ for certain permutations $\tau$ in $S_m$, namely
for a set of permutations chosen from the right cosets of the
subgroup $Gal(p')$.


Choose $\tau$ so that
$t - \tau(\beta ')$ divides $P_*$.
This means that
$t - \tau(\beta ') = t - \sigma(\beta ')$
for some permutation in $Gal(p')$ and hence
that $\tau \circ \sigma^{-1}$ is a permutation that permutes
the indices of repeated roots.
We may now reorder the repeated roots to define a new embedding
of $Gal(p')$ so that  $\tau \circ \sigma^{-1} = 1.$
Thus $t - \beta'$ is a factor of $P_*$ and it follows as before
that $P_*$ divides $Q_1'$ and that $T$ is a subgroup of $S$, as desired.

The last assertion of the theorem follows from the fact that a
subgroup of a soluble group is soluble. (See \cite{ste}.)
\end{proof}

The non-monic case of the last theorem can be deduced with the
following change of variables argument.

Suppose that $p$ is an irreducible polynomial in
$\mathbb{Q}[d][t]$ with non-zero specialisation $p'$. Choose a
rational number $a$ so that $p'(a) \ne 0$, and hence $p(a) \ne 0$.
Define the irreducible polynomial
$$
q(z) = t^np(t^{-1} + a)\frac{1}{p(a)}.
$$
Then $q$ is monic with well-defined specialisation
$$
q'(t) = t^np'(t^{-1}+a)\frac{1}{p'(a)}.
$$
The splitting fields of $p$ and $q$ are isomorphic as are those of
$p'$ and $q'$ and so it follows from the theorem above that
$Gal(p')$ is a subgroup of $Gal(p)$.

It is clear that the arguments above extend verbatim to the
specialisation of algebraic independents over any field of
characteristic zero and we shall need results in this setting. Let
$\mathbb{E}$ be such a field and let $\{d\}$ be a set of
algebraically independent variables over $\mathbb{E}$ with
rational field extension $\mathbb{E}(d)$.

\begin{thm}\label{T:galois2}
Let $p \in   {\mathbb E}  [{d}]  [t]$ be an irreducible polynomial
with Galois group ~~Gal$(p)$ when viewed as a polynomial in
${\mathbb E}   ({d})  [t]$. Let ${d}' \in {\mathbb E}^n$ be a
specialisation of ~${d}$ and let ~$p'$ be the associated
specialisation of $p$ with Galois group Gal$(p')$ over ${\mathbb
E}$. If $p'$ is non-constant  then Gal$(p')$ is a subgroup of
Gal$(p)$. In particular if $p$ is a radical polynomial then so too
is $p'$.
\end{thm}
\section{Planar $3$-connected CAD graphs are non-soluble}

We are now able to prove the main theorem stated in the
introduction.

Suppose, by way of contradiction, that there exists a maximally
independent $3$-connected planar graph which is soluble. Let $G$
be such a graph with the fewest number of vertices. We show that
$G$ is the doublet graph and that the doublet graph is not
soluble by radicals. This contradiction completes the proof.

By the reduction step, Theorem \ref{T:reduction}, the vertex minimal
graph $G$ has no
edge contraction to a $3$-connected maximally independent
planar graph. It thus follows from  the main reduction  theorem for
such graphs, Theorem \ref{T:4.1}, that either $|G| = 6$, and $G$
is the doublet (since $G$ is planar), or that $G$ has a proper
vertex induced maximally independent subgraph with three vertices
of attachment. However    minimality
rules out the latter possibility because  the next 
proposition  shows that such a proper subgraph
admits substitution by a smaller graph, namely a triangle,
and the resulting graph is soluble if $G$ is soluble.

\begin{prop}\label{P:8.1}
 Let $G$ be a 3-connected, maximally independent graph and let H be
a maximally independent subgraph of $G$ with $3$ vertices of
attachment $v_1, v_2$  and $v_3$. Let $G'$ be the graph which is
obtained from $G$ by deleting all the internal vertices of $H$ and
all the edges of $H$ and adding the edges $(v_1v_2), (v_2v_3),
(v_3v_1).$ Then $G'$ has the properties:

(i) $G'$ is $3$-connected.

(ii)    $G'$ is maximally independent.

(iii) If  the dimensions in the constraint 
equations defined by $G$ are
chosen as algebraic independents then the dimensions in the
equations defined by $G'$ are also algebraic independents.
\end{prop}

\begin{proof}
If $|H| = 3 $ then $H$ is the $3$-cycle and $G = G'$, so assume
$|H| \ge 4$. Note that  $H$ is connected since otherwise $G$ is not 
even be $2$-connected.

Every path in $G'$ derives from a path in $G\backslash H $ plus  paths
in $H$ which replaces segments $v_i \to v_j$ or $v_i \to v_j
\to v_k$ for $v_i, v_j, v_k$ chosen from the vertices of
attachment. For any set of independent paths in $G'$, at most one
of them contains any of the edges $(v_1v_2), (v_2v_3)$ or
$(v_3v_1)$. Thus every set of independent paths in $G'$ gives a
set of independent paths in $G$ and  (i) follows.

If $H_1$ and $H_2$ are any two edge disjoint subgraphs in $G$ then it follows easily
from the definition of $free(H)$ that
\[
free(H_1 \cup H_2) = free(H_1) + free(H_2) + 3 - 2|H_1 \cap H_2|.
\]
This gives immediately that $free(G')=0$. If $G'$ is not independent then there
is a subgraph $R$ of $G'$ with $free(R) < 0$ and there is an edge $(v_1v_2)$, 
say, which
is in $R$ but not in $G$. If $v_3$ is not in $R$ then 
$(R\backslash (v_1v_2)) \cup H$ is in $G$ 
and
$free((R\backslash(v_1v_2)) \cup H) < 0$ which contradicts 
the independence of $G$. If $v_3$ is in
$R$ then $(R\{(v_1v_2),(v_2v_3),(v_1v_3)\}) \cup H$ is in $G$ and
$free((R\{(v_1v_2),(v_2v_3),(v_1v_3)\}) \cup H) < 0$ 
which contradicts the independence of
$G$.

Theorem 2.3 implies that for algebracially independent dimensions $\{d_i\}$, any
zero of the variety of $G$ has coordinates $\{x_j\}$ which are algebraically
independent. This zero of the variety of $G$ gives a zero of the variety of $G'$
(with the same $\{x_j\}$ where they occur and with the same $\{d_i\}$ where they occur
and $d_{12}, d_{23}$ and $d_{13}$ computed from 
$d_{ij}=(x_i-x_j)^2+(y_i-y_j)^2)$ and this zero
therefore has coordinates which are algebraically independent. It follows
from Theorem 2.3 that the dimensions of $G'$ are algebraically independent.

\end{proof}

We now show that the doublet is a non-soluble CAD graph.\\

Let $v_{1} = (0,0), v_{2} = (1,0)$ be the vertices of the base
edge.  Introduce the coordinates $(x_{i}, y_{i})$ for the
remaining vertices $v_{i}, 3 \leq i \leq 6$, and the  dimensions
$d_{j}, 2 \leq j \leq 9$, for the non-base edges.  The indexing
scheme is illustrated in Figure 9.

\medskip
\begin{figure}
\centering
\includegraphics[width=9cm]{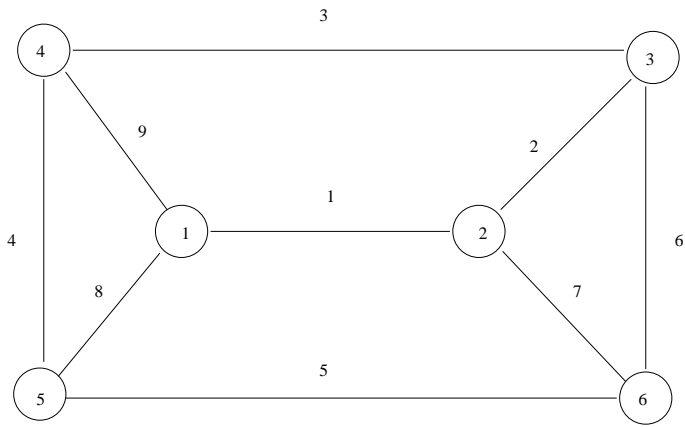}
\caption{Coordinatisation of the doublet.}
\end{figure}
\medskip

The resulting polynomials $\{f \}$ for the normalised constraint
equations take the form
\medskip
$$\begin{array}{ccccc}
\vspace{2mm}
x_{4}^2 &  +  & y_{4}^{2}  & - & d_{9}^{2}\\
\vspace{2mm}
x_{5}^{2} &  +  & y_{5}^{2}  & - & d_{8}^{2}\\
\vspace{2mm}
(x_{3} - 1)^2 &  +  & y_{3}^{2} & - & d_{2}^{2}\\
\vspace{2mm}
(x_{6} - 1)^2 &  +  & y_{6}^{2} & - & d_{7}^{2}\\
\vspace{2mm}
(x_{3} - x_{4})^2 &  +  & (y_{3} - y_{4})^{2} & - & d_{3}^{2}\\
\vspace{2mm}
(x_{4} - x_{5})^2 &  +  & (y_{4} - y_{5})^{2} & - & d_{4}^{2}\\
\vspace{2mm}
(x_{5} - x_{6})^2 &  +  & (y_{5} - y_{6})^{2} & - & d_{5}^{2}\\
\vspace{2mm}
(x_{6} - x_{3})^2 &  +  & (y_{6} - y_{3})^{2} & - & d_{6}^{2}.
\end{array}$$
\medskip

For each   choice of real algebraically independent squared
dimensions ${\displaystyle d_{2}^{2}, \ldots, d_{9}^{2}}$ these
equations determine a zero-dimensional complex affine variety
$V(\{ f \} )$ in ${\mathbb C}^8$.

Note that the fifth equation, and its three successors, admit the squared form
$$(d_{3}^{2} -  (x_{3} - x_{4})^{2}  +  y_{3}^{2}  +
y_{4}^{2})^{2} - 4y_{3}^{2} \, y_{4}^{2} = 0,$$ which in turn yields an
equation in $x_{3}$ and $x_{4}$ alone on substituting for
$y_{3}^{2}$ and $y_{4}^{2}$ from the first four equations.  In
this way we obtain a system $\{g \} = \{ g_1 , g_2 , g_3, g_4 \}$
of four quartic equations in $x_3, x_4, x_5, x_6$ and the squared
dimensions.  It follows that the projection $\pi \left(V( \{f
  \} \right)$ for the variables $x_3, x_4, x_5, x_6$ is a subset of
the variety $V \left( \{ g \} \right)$ in $\mathbb{C}^4$.

To see that the doublet graph is (generically) non soluble we show
first there is a specialised integral dimensioned doublet which
has non radical solutions.  This is achieved by a Maple
calculation of successive resultants of the associated specialised
constraint equations $\{ g' \}$;

$$\begin{array}{lll}
\vspace{1mm}
h_1' & = & Res (g_1', g_2', x_4),\\
\vspace{1mm}
h_2' & = & Res (g_3', g_2', x_6),\\
\vspace{1mm}
h_3' & = & Res (h_1', h_2', x_5).
\end{array}$$

This results in an integral univariate polynomial $h_3' (x_3)$
which lies in the ideals $I( \{ f' \} )$ and $I( \{ g' \} )$. The
polynomial $h_3'$ is of degree $28$ which normally rules out
convenient computer algebra calculation of the Galois group.
However for our well-chosen dimension values (determined by judicious trial
and error) the polynomial factors as  a product of four
irreducible polynomials of degrees $6$, $6$, $8$, $8$. The Galois
groups of these polynomial factors are computed in the Appendix, and 
each is a
full symmetric group. It follows that $h_3'$ and $V
\left( \{ f' \} \right)$ are
not radical over ${\mathbb Q}$.

\begin{thm}\label{T:Zdoublet}
There exists an integral dimensioned doublet graph which is not soluble by radicals.
\end{thm}

\begin{proof}
With the labelling order above consider the unsquared dimensions
$1, 5, 15, 10, 16, 8, 5, 13, 13$. (The two triangles in this integral doublet
are isosceles, with sides $10, 13, 13$ and $8, 5, 5$.) By the Appendix 
$h_3'$ is a non-radical polynomial.
\end{proof}

We now use the Galois group specialisation theorem to show that
the doublet graph is generically non-soluble. The generic
polynomial $h_3$ is not conveniently computable but we examine the
resultant calculation more closely to see that $h_3'$ is the
specialisation of the corresponding resultant polynomial $h_3$ for the generic
equation set.\\

\begin{lma} Let $f_1, f_2$ be polynomials in $ \{ x \}, \{ d \}$ viewed
as polynomials in $ \{ x \}$ with coefficients in ${\mathbb E} (
\{ d \} )$.  Let $\{ d' \}$ be a specialisation resulting in
specialisations $f_1', f_2'$ such that $\deg (f_{i}, x_1) = \deg
(f_i', x_1)$ for $i=1,2$.  Then the specialisation of $ \mbox{Res}
(f_1, f_2, x)$ is equal
to $ \mbox{Res} (f_1', f_2', x)$.\\
\end{lma}

\begin{proof} Immediate on examination of the definition of the resultant
as a Sylvester determinant.
\end{proof}

For our polynomial equations $\{g\}$ a simple Maple verification
shows that if $h_1 = \mbox{Res} (g_1, g_2, x_4), h_2 = \mbox{Res}
(g _3, g_4, x_6)$ then
$$ \deg(h_1,x_4)=\deg(h_1',x_4)=\deg(h_2,x_6)=\deg(h_2',x_6)=4.
$$

Although the polynomial $h_3$ is not readily computable the lemma
shows that $h_3'$ is the specialisation of $h_3$.

\begin{thm}
The doublet graph is non-soluble.\\
\end{thm}

\begin{proof} By Theorem \ref{T:Zdoublet} and its proof $h_3'$ is a non radical
polynomial and in fact all the zeros of its irreducible factors
are non radical over ${\mathbb Q}$.  By the Galois  group
specialisation theorem it follows that $h_3$ must be non radical
over ${\mathbb Q} ( \{d \})$ and the theorem follows.
\end{proof}

\newpage

{\bf Appendix}
\medskip

The polynomial $h_3'$ and its factors are computed by the
following Maple code.
\begin{verbatim}

 d2:= 13; d3:=  15; d4:= 8; d5:= 16;
 d6:= 10; d7:=  13; d8:= 5; d9:= 5;
 yy4:=d9^2-x4^2; yy5:=d8^2-x5^2;
 yy3:=d2^2-(x3-1)^2; yy6:=d7^2-(x6-1)^2;
 A:= (d3^2- (x3^2+x4^2 - 2*x3*x4 +  yy3 + yy4))^2 -4*yy3*yy4;
 B:= (d4^2- (x4^2+x5^2 - 2*x4*x5 +  yy4 + yy5) )^2-4*yy4*yy5;
 C:= (d5^2- (x5^2+x6^2 - 2*x5*x6 +  yy5 + yy6) )^2-4*yy5*yy6;
 E:= (d6^2- (x6^2+x3^2 - 2*x6*x3 +  yy6 + yy3) )^2-4*yy6*yy3;
 eqns:={A=0,B=0,C=0,E=0}; expand(eqns);
 X:=resultant(A,B,x4): Y:=resultant(C,E,x6):
 Z:=resultant(X,Y,x5):
 factor(Z):

\end{verbatim}
The irreducible factors are the following four integral
polynomials and (according to Maple) each is non-soluble over
$\mathbb{Q}.$

\begin{small}
$$
\begin{array}{l}
731161600000x_3^8 - 2884724544000x_3^7 - 254604702168560x_3^6 +
\\
929745074065696x_3^5 + 29180343859430360x_3^4 - 104245652941659832x_3^3 -
\\
1119855862049129679x_3^2 + 4022769219537416744x_3
 + 1620713038685642896,
\end{array}
$$

$$
\begin{array}{l}
731161600000x_3^8 - 5275493184000x_3^7 - 202247115019760x_3^6 + \\
1002422141698336x_3^5 + 16575444136627160x_3^4 - 46366435207277752x_3^3 - \\
299095702632348879x_3^2 + 813935120915198504x_3 +
13663404945744016,
\end{array}
$$

$$
\begin{array}{l}
753831936x_3^6 - 84641660928x_3^5 - 4996031627504x_3^4
+ \\
486105086115256x_3^3 +  36795384322988721x_3^2 +
920226256962743080x_3 + \\
 10127898920872530064,
\end{array}
$$

$$
\begin{array}{l}
2747437056x_3^6 + 143122194432x_3^5 - 17613405584624x_3^4 -
615688594921544x_3^3 +
\\
69050497529701041x_3^2 - 776224290995754200x_3 +
1152246393155768464.
\end{array}
$$
\end{small}

\end{document}